\documentclass{article}

\pdfoutput=1

\usepackage[utf8]{inputenc}

\usepackage[a4paper, total={6in, 9in}]{geometry}
\usepackage{authblk}

\usepackage{amssymb,amsfonts,amsmath}
\usepackage{subfigure}
\usepackage{url}
\usepackage{xcolor}
\usepackage[linesnumbered,ruled]{algorithm2e}

\usepackage{scalefnt}

\usepackage{hyperref}

\usepackage{amsthm}

\usepackage{nccmath}

\usepackage{multirow}
\usepackage{colortbl}
\usepackage{booktabs}

\usepackage{adjustbox}
\usepackage{setspace}

\usepackage{tikz}
\usepackage{pifont}
\usetikzlibrary{calc, fit, math, decorations, backgrounds}

\usepackage{subfigure}

\newtheorem{rules}{Rule}

\newtheorem{remark}{Remark}

\usepackage{acro}
\acsetup{
  list/sort = true ,
  make-links =true
}
\DeclareAcronym{SSB}{
  short = SSB,
  long  = Static Symmetry Breaking ,
  tag = abbrev
}

\DeclareAcronym{BPM}{
  short = BPM,
  long  = Batch Processing Machines ,
  tag = abbrev
}
\DeclareAcronym{MILP}{
  short = MILP ,
  long  = Mixed Integer Linear Programming ,
  tag = abbrev
}
\DeclareAcronym{BeB}{
  short = B\&B,
  long  = Branch-and-Bound ,
  tag = abbrev
}

\title{Arc-flow approach for single batch-processing machine scheduling}

\author[1]{Renan Spencer Trindade\thanks{Corresponding author.\\ 
Address: 1 rue Honor\'e d'Estienne d'Orves, 91120  Palaiseau, France. Email: rst@lix.polytechnique.fr}}
\author[2]{Olinto C. B. de Araújo}
\author[3]{Marcia Fampa}

\affil[1]{LIX, CNRS, \'Ecole Polytechnique, Institut Polytechnique de Paris}
\affil[2]{Colégio Técnico Industrial de Santa Maria, Universidade Federal de Santa Maria}
\affil[3]{Instituto de Matem\'atica, COPPE, Universidade Federal do Rio de Janeiro}

\date{}
\usepackage{graphicx}

\begin{document}

\maketitle

\begin{abstract}
    We address the problem of scheduling jobs with non-identical sizes and distinct processing times  on a single batch processing machine, aiming at minimizing the makespan. The extensive literature on this  NP-hard problem mostly focuses on heuristics. Using an arc flow-based optimization approach, we construct an ingenious  formulation that represents it as a  problem of determining flows in graphs. The size of the formulation increases with the number of distinct sizes and processing times among the jobs, but it does not increase with the number of jobs, which makes it very effective to solve large instances  to optimality, especially when multiple jobs have equal  size and processing time. We compare our model to other models from the literature showing its clear superiority on benchmark instances, and proving optimality of random  instances with up to 100 million~jobs. 
\end{abstract}


\section{Introduction}
This paper addresses  \ac{BPM} scheduling, which has been extensively explored in the literature, motivated by applications in industries and by the challenging solution of real-world problems. The main goal in \ac{BPM}  scheduling  problems is to group jobs into batches and process them simultaneously in a  machine, to facilitate the tasks and to reduce the time spent in handling the material. Although there are many variations of problems involving \ac{BPM}, the version considered in this work is particularly suitable to model scheduling problems that arise in reliability tests in the semiconductor industry, in operations called burn-in \cite{UZSOY1994}.

The burn-in operation is used to test electronic circuits and consists of designating them to industrial ovens, submitting them to thermal stress for a long period. The test of each circuit is considered here as a job and requires a minimum time inside the oven, which is the processing machine.  The jobs cannot be processed directly in the machine, they need to be placed on a tray, respecting its capacity. Each group of jobs assigned to a tray is considered a batch.  The minimum time of the circuit inside the oven is set a priori, based on the supplier requirements. It is assumed that it is possible to keep the circuit in the oven longer than necessary, with no prejudice, but the circuit cannot be removed from the oven before its required processing time is fulfilled. Therefore, the processing time of a batch is determined by the longest processing time among all jobs assigned to it. The burn-in tests are a bottleneck in final testing operations, and the efficient scheduling of these operations aims at maximizing productivity and reducing flow time in the stocks, which are  major concerns for management. In \cite{Lee1992}, \cite{Tai2008}, and \cite{Chung2009},
the complexity of burn-in tests on electronic circuits and liquid crystal displays is discussed, 
reinforcing the importance of an efficient scheduling for them.
  
Compared to the history of the semiconductor manufacturing,  the interest of  researchers in \ac{BPM} is recent.   A  review  on scheduling models considering batch processing machines is presented in  \cite{Potts2000}, where the authors point out that the research effort in designing algorithms for these problems is worthwhile.
A survey on \ac{BPM} problems is presented  in \cite{Mathirajan2006}, considering publications between 1986 and 2004 (part of 2004 only). A more recent  survey is published in \cite{Monch2011}. 
  
The \ac{BPM} problem  specifically addressed in this paper is  the single batch processing machine scheduling problem, denoted by $1|s_j,B|C_{\max}$, where jobs have non-identical sizes and processing times, and the objective is to minimize the makespan. The problem was reported for the first time in \cite{UZSOY1994}, where its NP-hard complexity is proven and a heuristic approach to solve it is proposed. Heuristics are also proposed for this problem in \cite{Ghazvini1998}, where instances with up to 100 jobs are considered.  Two approximation algorithms are presented in \cite {Zhang2001} with approximation ratios of $3/2$ and $7/4$ of the optimal solution, in the worst case.
In \cite{Melouk2004}, the simulated annealing meta-heuristic was applied to $1|s_j,B|C_{\max}$ and a \ac{MILP} formulation was presented for the problem.  This work also proposes configurations for test instances that were widely used in later works. Computational results are shown for instances with up to 100 jobs,  comparing the heuristic solutions to the solutions obtained with the \ac{MILP} formulation. 
Other meta-heuristics are also applied to problem $1|s_j,B|C_{\max}$ in the literature, namely, genetic algorithm  (\cite{Damodaran2006} and \cite{Kashan2006}), tabu search (\cite{Meng2010}),  and GRASP (\cite{Damodaran2013}). These four papers  consider instances with up to 100 jobs as well.
In addition, the bee colony meta-heuristic is also applied to the problem in \cite{Al-Salamah2015}, where results for instances with up to 200 jobs are shown.
In \cite{Chen2011}, a heuristic based on a special case of the clustering problem is proposed, and test instances with up to 500 jobs are considered.
In \cite{Li2017}, an enumeration scheme for heuristics is proposed. The work uses First Fit Longest Processing Time (FFLPT) and Best Fit Longest Processing Time (BFLPT) considering identical job sizes, and instances with up to 500 jobs.
In \cite{Lee2013}, two heuristics are proposed based on a decomposition of the original problem, where relaxations of the problem are solved. Instances with up to 100 jobs are considered in this work. 

Most publications addressing problem  $1|s_j,B|C_{\max}$, apply heuristic and meta-heuristic approaches.  Nevertheless, exact methods are also applied in the literature, as for example, in \cite{RafieeParsa2010}, \cite{Trindade2018}, \cite{TeseTrindade}, and \cite{Muter2020}.  
In \cite{RafieeParsa2010},  a formulation for $1|s_j,B|C_{\max}$  is presented, using a partition problem in the context of Dantzig-Wolfe decomposition. A branch-and-price algorithm is applied to solve the problem to optimality. Instances with up to 500 jobs are considered in the computational experiments.
In \cite{Trindade2018} and \cite{TeseTrindade}, the authors consider the exact solution of the problem, improving a formulation from the literature. Sets of symmetric solutions for the formulation are identified. Symmetry breaking constraints are then proposed, taking into account specific properties of the problems and their optimal solutions to propose new stronger formulations and avoid undesirable symmetric solutions in their feasible sets.
\cite{Muter2020} proposes a column-and-cut generation (CCG) algorithm for the instances that could not be solved to optimality with the model proposed by \cite{Trindade2018}. The method consists of running the model proposed by \cite{Trindade2018} for 300 seconds to obtain an upper bound or, if possible, an optimal solution, and then running the GCC method to improve the bounds.
This work also proposes an exact algorithm for the version of the problem that considers parallel machines.

In this paper, we aim at solving problem $1|s_j,B|C_{\max}$  to optimality as well, but considering  an arc-flow optimization approach.  
Arc flow-based  approaches  have been applied to  classical combinatorial optimization problems and have allowed formulations with a pseudo-polynomial number of variables and constraints. For example, see \cite{DeCarvalho1998}, where the author proposes a branch-and-price algorithm for an arc flow-based formulation for a cutting-stock problem. Later, the methodology was extended for the bin-packing problem in \cite{ValeriodeCarvalho1999}.
An alternative arc-flow formulation for the cutting-stock problem is proposed in \cite{Brandao2016} and \cite{Brandao2017}, which uses graph compression technique.
These formulations were recently tested and compared in \cite{Delorme2016} against several other models and problem-specific algorithms on one-dimensional bin packing and cutting stock problems. The results show that the arc flow-based formulations outperform all other models.
\cite{MACEDO2010991} extended 
the model proposed by Valério de Carvalho for the Two-Dimensional cutting stock problem. The authors tested the resulting model in a set of real instances of the furniture industry, and it proved to be more efficient than other previous methods
%

%
In   \cite{Martinovic2018} the arc-flow model and the one-cut model are compared for the one-dimensional cutting-stock problem, and reduction techniques for both approaches are presented. 
Regarding scheduling problems, we are aware of only two papers that apply an arc-flow approach. In \cite{Kramer2018}  problem $P||\sum W_jC_j$ is considered, where the goal is to schedule a set of jobs on a set of identical parallel machines, while minimizing the total weighted completion time. In \cite{Mrad2018}, the makespan minimization problem on identical parallel machines, $P||C_{max}$, is considered. Differently from what we present, these works do not consider important features  of scheduling problems, such as batch processing machines, non-identical job sizes, and machine capacity,  they consider different scheduling problems.

Motivated by the good results obtained with arc-flow approaches for combinatorial problems,  we  present in this paper an ingenious arc flow-based formulation for problem $1|s_j,B|C_{\max}$, which represents it as a problem of determining flows in graphs. A special graph $G$, and  arc-flow structures corresponding to the distinct processing times of the jobs, are introduced to form the core of our formulation. With the approach used in the modeling, we are able to propose a formulation, where neither the number of variables  nor of constraints increases with the number of jobs, they increase only with the number of distinct sizes and processing times among the jobs. 
Another important feature of the formulation is that it does not contain in its feasible set, some of the symmetries that are contained in the feasible sets of the formulations presented in \cite{RafieeParsa2010} and \cite{Trindade2018}.
We present numerical experiments, where very large instances of the problem, with a number of jobs never considered in the literature,  are solved to optimality in short time, even if compared to the solution time of much smaller instances, when using  other formulations from the literature. After the development of this work, we could also successfully extend the idea of the arc-flow model for the case of parallel machines.  Computational results for this extension can be found in  \cite{ISCO2020}.

We organize the paper as follows:
In Section \ref{chap:3}, we introduce problem $1|s_j,B|C_{\max}$  and present two formulations from the literature.
In Section \ref{ArcFlow}, we present an arc flow-based formulation for the problem.
In Section \ref{ComputationalResults}, we discuss our numerical experiments comparing the arc-flow formulation to formulations from the literature. 
In Section \ref{sec:conc}, we present some concluding remarks and discuss future work. 

\section{Problem definition} 
\label{chap:3}
The single batch processing machine scheduling problem, or problem  $1|s_j,B|C_{\max}$, can be formally defined as follows. Given a set $J:=\{1,\ldots,n_J\}$ of jobs, each job $j \in J$ has a processing time $p_j$ and a size $s_j$. Each job must be assigned to a batch $k \in K:=\{1,\ldots,n_K\}$, respecting the given capacity  $B$ of the machine, so that the sum of the sizes of the jobs assigned to a batch cannot exceed $B$. We assume that $s_j\leq B$, with   $s_j,B\in\mathbb{Z}$, for all $j\in J$. The batches must be all processed in a single machine, one at a time, and all the jobs assigned to a single batch are processed simultaneously. The processing time $C_k$ of each batch $k \in K$ is defined as longest processing time among all jobs assigned to it, i.e., $C_k := \max \{ p_j: j \mbox{ is assigned to } k \} $. Jobs cannot be split between batches. It is also not possible to add or remove jobs from the machine while the batches are being processed. The goal is to schedule the jobs on the batches so that the makespan ($C_{\max}$) is minimized, where  the makespan is defined as the time required to finish processing the last batch.
  
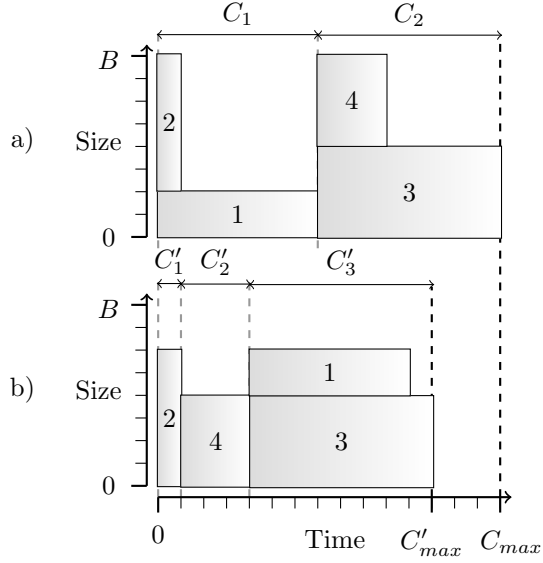
\begin{figure}[!ht]
  \def\scale{0.3}
  \begin{center}
    \newcommand{\GenerateXAxis}[0]{
    \draw [->, thick] (0,0-0.5) node [yshift = -0.5cm] {0} -- node [yshift = -0.8cm] {\NameXAxis} (\TimeMaxValue+0.5,0-0.5) node [yshift = -0.6cm] {$C_{max}$};
    \foreach \i in {1,...,\TimeMaxValue}
    {
        \draw[-] (\i,-0.5) -- (\i,-1);
    }
    \draw[-, thick] (0,-0.5) -- (0,-1.18);
    \draw[-, thick] (\TimeMaxValue,-0.5) -- (\TimeMaxValue,-1.18);
}

\newcommand{\GenerateSizeAxe}[1]{
    \draw [->, thick] (0-0.5,\SizeBeginValueY + 0) node [xshift = -0.5cm] {0} -- node [xshift = -0.2cm] (#1) {Size} (0-0.5,\SizeBeginValueY + \SizeMaxValue+0.5) node [xshift = -0.5cm, yshift = -0.2cm] {$B$};
    
    \foreach \i in {1,...,\SizeMaxValue}
    {
        \draw[-] (-0.5,\SizeBeginValueY + \i) -- (-1,\SizeBeginValueY + \i);
    }
    \draw[-, thick] (-0.5,\SizeBeginValueY) -- (-1.18,\SizeBeginValueY);
    \draw[-, thick] (-0.5,\SizeBeginValueY+\SizeMaxValue) -- (-1.18,\SizeBeginValueY+\SizeMaxValue);
}

\newcommand{\CreateNewJob}[4]{
    \node[fit={#4 ($#4 +
    #1
    $)}, inner sep=0.3pt, draw=black, shade, left color=mycolor, right color=white, label=center:#3] (#2) {};
}

\begin{tikzpicture}[scale=\scale, shorten >=0pt,node distance=0cm,auto]
\definecolor{mycolor}{RGB}{220,220,220}
\tikzmath{
    \SizeMaxValue = 8;
    \TopLineBegin = \SizeMaxValue+1;
    \SizeBeginValueY = 0;
    \DistanceBtweenGraphs = 11;
    \BatchRelease = 0;
}

\draw[-, thick, dashed] (15,-0.5) -- (15,8.8+\DistanceBtweenGraphs*1);
\draw[-, thick, dashed] (12,-0.5) -- (12,8.8+\DistanceBtweenGraphs*0);

\draw[-, thick, dashed, opacity=0.4] (0,-0.5+\DistanceBtweenGraphs) -- (0,\TopLineBegin+\DistanceBtweenGraphs);
\draw[-, thick, dashed, opacity=0.4] (7,-0.5+\DistanceBtweenGraphs) -- (7,\TopLineBegin+\DistanceBtweenGraphs);

\draw[-, thick, dashed, opacity=0.4] (0,-0.5+\DistanceBtweenGraphs*0) -- (0,\TopLineBegin+\DistanceBtweenGraphs*0);
\draw[-, thick, dashed, opacity=0.4] (1,-0.5+\DistanceBtweenGraphs*0) -- (1,\TopLineBegin+\DistanceBtweenGraphs*0);
\draw[-, thick, dashed, opacity=0.4] (4,-0.5+\DistanceBtweenGraphs*0) -- (4,\TopLineBegin+\DistanceBtweenGraphs*0);

\tikzmath{
    \TimeMaxValue = 15;
    \NameXAxis = "Time";
}
\GenerateXAxis
\draw[-, thick] (12,-0.5) node [yshift = -0.6cm] {$C_{max}'$} -- (12,-1.18);

\GenerateSizeAxe{batch1}

\tikzmath{
\SizeBeginValueY = \DistanceBtweenGraphs;
}
\GenerateSizeAxe{batch2}

\node [align=left, left of=batch1, xshift=-1cm] {b)};
\node [align=left, left of=batch2, xshift=-1cm] {a)};

\tikzmath{
    \SizeBeginValueY = \DistanceBtweenGraphs;
}
\CreateNewJob{(7,2)}{job_1}{1}{(0,\SizeBeginValueY)};
\CreateNewJob{(1,6)}{job_2}{2}{($(job_1.north west) + (1pt,0)$)};
\draw[<->] ($(job_1.south west) + (0,\TopLineBegin)$) -- node [above] {$C_1$} ($(job_1.south east) + (0,\TopLineBegin)$);

\CreateNewJob{(8,4)}{job_3}{3}{($(job_1.south east) + (-1pt,1pt)$)};
\CreateNewJob{(3,4)}{job_4}{4}{($(job_3.north west) + (1pt,0)$)};
\draw[<->] ($(job_3.south west) + (0,\TopLineBegin)$) -- node [above] {$C_2$} ($(job_3.south east) + (0,\TopLineBegin)$);

\tikzmath{
    \SizeBeginValueY = 0;
}
\CreateNewJob{(1,6)}{job_22}{2}{(0,\SizeBeginValueY)};
\draw[<->] ($(job_22.south west) + (0,\TopLineBegin)$) -- node [above] {$C_1'$} ($(job_22.south east) + (0,\TopLineBegin)$);

\CreateNewJob{(3,4)}{job_42}{4}{($(job_22.south east) + (-1pt,1pt)$)};
\draw[<->] ($(job_42.south west) + (0,\TopLineBegin)$) -- node [above] {$C_2'$} ($(job_42.south east) + (0,\TopLineBegin)$);

\CreateNewJob{(8,4)}{job_32}{3}{($(job_42.south east) + (-1pt,1pt)$)};
\CreateNewJob{(7,2)}{job_12}{1}{($(job_32.north west) + (1pt,0)$)};
\draw[<->] ($(job_32.south west) + (0,\TopLineBegin)$) -- node [above] {$C_3'$} ($(job_32.south east) + (0,\TopLineBegin)$);

\end{tikzpicture}
  \caption{Feasible solutions with $(a)$ $2$ batches  and $(b)$ $3$ batches.} 
  \label{fig_contraexemplo}
  \end{center}
\end{figure}

The  number of batches used in the solution of $1|s_j,B|C_{\max}$ is not fixed and should be optimized. It will depend on the number of jobs, their sizes,  and the machine capacity. In this section, we present two formulations for the problem from the literature, where $n_K$ actually represents the maximum number of batches allowed in the solution.   We will see that  this number interferes in the problem size.  
  In the worst case, each job is assigned to a different batch, and the number of batches used is equal to the number of jobs $n_J$.
It is important to note that even when  it is possible to find a feasible solution for a problem where the number of batches used is $\tau < n_J$, we should not limit the maximum number of batches to $\tau$ in the problem's formulation, as we should  consider the case in which creating an additional batch in the solution would decrease the makespan.
As an example of this situation, we show in Figure \ref{fig_contraexemplo} two feasible solutions for an instance of  $1|s_j,B|C_{\max}$. Although the number of batches used in solution $(a)$ is smaller than in solution $(b)$, its makespan  is bigger.
  Therefore, to assure the correctness of the formulations presented next,  it is always considered $n_K=n_J$.

Let us consider the following decision variables for all $j\in J$, $k\in K$: 
\[
	\begin{array}{lll}
  x_{jk} &= \left\{ \begin{array}{ll}
	1, & \mbox{if job $j$ is assigned to batch $k$},\\
	0, & \mbox{otherwise},\end{array} \right.\\
  y_{k} &= \left\{ \begin{array}{ll}
  1, & \mbox{if batch $k$ is used},\\
  0, & \mbox{otherwise},\end{array} \right. \\
  C_k&\multicolumn{2}{l}{ \; :\;\; \mbox{processing time of batch  } k.}\\
  \end{array}
\]
		
In \cite{Melouk2004} the following MILP formulation  is proposed for problem $1|s_j,B|C_{\max}$. Other very similar formulations and sometimes this exactly same one, are used in other papers as a comparative basis in computational experiments.
 
\begin{fleqn}
\begin{align} \hspace{0.7cm}
  (\mbox{MILP}_1)\;\;&\min \; \sum_{k \in K} C_{k}, &                                                  \label{Model1_FO}\\
    &\sum_{k \in K} x_{jk} = 1 ,		    & \forall j \in J,                        \label{Model1_Assig} \\
    &\sum_{j \in J} s_j x_{jk} \le B y_{k},	& \forall k \in K,                        \label{Model1_Capac} \\
    & C_k \geq p_j x_{jk}, 				    & \forall j \in J, \forall k \in K,       \label{Model1_ProcTim} \\
    & x_{jk} \leq y_k, 				    & \forall j \in J, \forall k \in K,           \label{Model1_Red} \\
    & y_k \in \{0,1\},  						& \forall k \in K,                            \label{Model1_DomYk} \\
    & x_{jk} \in \{0,1\}, 				& \forall j \in J, \forall k \in K.        \label{Model1_DomXjk}             
  \end{align}
\end{fleqn}
  
The objective function \eqref{Model1_FO} minimizes the makespan, given by the sum of the processing times of all batches. 
Constraints  \eqref{Model1_Assig} determine that each job is assigned to a single batch. Constraints \eqref{Model1_Capac} determine that the machine capacity is respected. 
Constraints \eqref{Model1_ProcTim} determine the processing times of the batches.
Note that constraints \eqref{Model1_Red} are redundant because of \eqref{Model1_Capac}, but are added to strengthen the linear relaxation of the formulation.

In \cite{Trindade2018}, the authors address the symmetry in formulation  (MILP$_1$). They present a new formulation for the problem, where symmetric solutions are eliminated from the feasible set of  (MILP$_1$),  with the following approach. Firstly,  the indexes of the jobs are defined by ordering them by their processing times. More specifically, it is considered that $p_1 \leq p_2 \leq \ldots \leq p_{n_J}$. Secondly,  it is determined that batch $k$ can  only be used if job $k$ is assigned to it, for all $k\in K$. Thirdly, it is determined that job $j$ can only be assigned to batch $k$ if $j\leq k$.
Considering the above,  the following formulation for $1|s_j,B|C_{\max}$ is proposed in \cite{Trindade2018}:

\begin{fleqn}
    \begin{align} \hspace{0.7cm}
     (\mbox{MILP}_1^+)\;\;&\min \; \sum_{k \in K} p_{k}x_{kk},    {\;\;\;\;\;\;\;\;\;\;\;\;\;\;\;\;\;\;\;\;\;\;\;\;}                                         \label{Model2_FO}\\
    &\sum_{k \in K  :  k\geq j} x_{jk} = 1, 		& \forall j \in J,                 \label{Model2_Assig2}\\
    &\sum_{j \in J :  j \leq k} s_j x_{jk} \le B x_{kk},	& \forall k \in K,     \label{Model2_Capac}\\
    & x_{jk} \leq x_{kk}, 				& \forall j \in J, \forall k \in K  :  j \leq k , \label{Model2_Red2}\\
    & x_{jk} \in \{0,1\},				& \forall j \in J, \forall k \in K  :  j \leq k.                 \label{Model2_DomX}
    \end{align}
\end{fleqn}
		 
The objective function \eqref{Model2_FO} minimizes the makespan, given by the sum of the processing times of the batches used. 
Constraints \eqref{Model2_Assig2} determine that each job $j$ is assigned to a single batch $k$, such that $k\geq j$. Constraints \eqref{Model2_Capac} determine that the machine capacity  is respected. They also ensure that each batch $k$ is used if and only if job $k$ is assigned to it. Constraints \eqref{Model2_Red2} are redundant together with \eqref{Model2_Capac}, but are included to strengthen the linear relaxation of the model. 

We note here that there are similarities between problem $1|s_j,B|C_{\max}$ and the classic bin-packing and cutting-stock combinatorial optimization problems. \cite{UZSOY1994} proves that the bin-packing problem is a special case of $1|s_j,B|C_{\max}$, in which the  processing times $p_j$, for all jobs,   are equal to one. Therefore, any formulation for $1|s_j, B|C_{\max}$ can also model the bin-packing problem. The opposite, on the other hand, is not valid, because formulations for bin packing cannot distinguish items with different processing times, neither calculate the makespan. Moreover, the cutting-stock problem can be seen as a generalization of the bin-packing problem, with unit demands, as shown in \cite{Vanderbeck1999} and \cite{Delorme2016a}. However,  it is not a generalization of problem $1|s_j, B|C_{\max}$. As the cost of each pattern in the cutting-stock problem is originally set to one, the different processing times of the batches in the  $1|s_j, B|C_{\max}$ problem, cannot be represented in the cutting-stock formulation.





\section{An arc flow-based optimization approach} \label{ArcFlow}

In this section we propose an arc flow-based optimization approach  to problem $1|s_j,B|C_{\max}$, formulating it  as a problem of determining flows in graphs. 
The arc-flow approach was also applied to formulate other problems, such as  the cutting-stock problem, in   \cite{DeCarvalho1998}, and  the bin-packing problem, in \cite{ValeriodeCarvalho1999}.
It is important to note, however, that the formulations for these problems cannot be trivially  applied to $1|s_j,B|C_{\max}$, as mentioned in the previous section. Some specific aspects of this last problem have to be exploited and considered in the new formulation that we propose in the following.

We initially define a directed  graph $G:=(V, A)$, in which each (integer) physical space of the batch with capacity $B$ is represented by a node, so  \begin{equation}\label{defV} V:=\{0,1,\ldots,B \}.\end{equation} 

The set of directed arcs  $ A $ is divided into three subsets: the set of \emph{job arcs} $ A^J$, the set of \emph{loss arcs} $A^L$, and the set with a \emph{feedback arc} $A^F$. Therefore, $A=A^J \cup A^L \cup A^F$. 

Each arc $(i,j)$ of the subset $A^J$ represents the existence of at least one job $k$ of size $s_k$, such that $s_k=j-i$. More specifically, the subset  $ A^J $ is  defined as:
\begin{align}
\label{arc_J}
A^J := \{(i,j):\exists\,  k \in J \ | \  s_k = j-i, \   i,j \in V , \ i<j\}.
\end{align}

As it will become clear with the presentation of our arc-flow formulation in Subsection \ref{sec:form}, to compose valid paths  and represent all possible solutions for the problem, it is necessary to include the \emph{loss arcs} in $G$, which represent empty spaces in the batches. More specifically, the subset of arcs $A^L$ is  defined as:
\begin{align}
\label{arc_L}
A^L := \{(i,B): i\in V\setminus\{0,B\} \}.
\end{align}

Finally, the \emph{feedback arc} is used to connect the last node to the first one, defined as:
\begin{align*}
A^F:=\{(B,0)\}.
\end{align*}

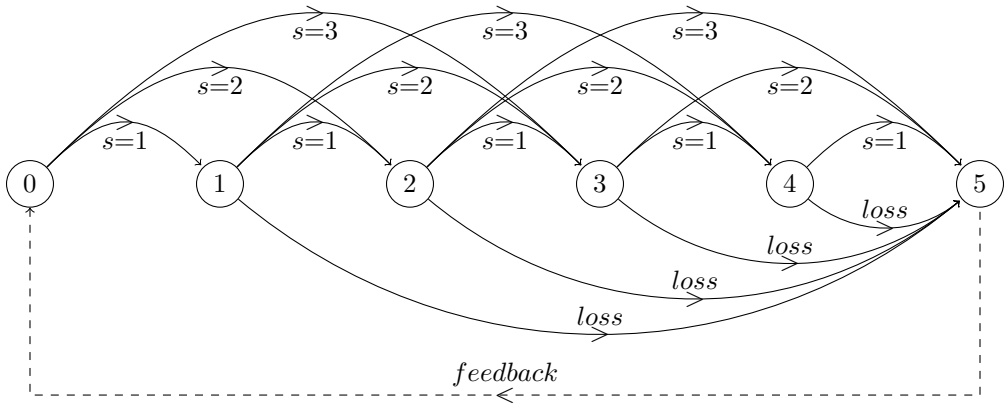
\begin{figure}[!ht]
	\centering
	\caption{Example of the graph $G$ representing the arc-flow structure used to model problem $1|s_j,B|C_{\max}$. }
	\label{fig_flow1}
\begin{tikzpicture}[shorten >=1pt,node distance=2.5cm,auto]
\node[circle, draw] (q_0) {$0$};
\node[circle, draw] (q_1) [right of=q_0] {$1$};
\node[circle, draw] (q_2) [right of=q_1] {$2$};
\node[circle, draw] (q_3) [right of=q_2] {$3$};
\node[circle, draw] (q_4) [right of=q_3] {$4$};
\node[circle, draw] (q_5) [right of=q_4] {$5$};
\draw[->] (q_0) to [looseness=1.4, bend left=45, 
                    edge node={node [sloped,below] {$s$=1}},
                    edge node={node [yshift = -0.23cm] {$>$}}](q_1);
\draw[->] (q_0) to [looseness=1.4, bend left=45, 
                    edge node={node [sloped,below] {$s$=2}},
                    edge node={node [yshift = -0.23cm] {$>$}}](q_2);
\draw[->] (q_0) to [looseness=1.4, bend left=45, 
                    edge node={node [sloped,below] {$s$=3}},
                    edge node={node [yshift = -0.23cm] {$>$}}](q_3);
\draw[->] (q_1) to [looseness=1.4, bend left=45, 
                    edge node={node [sloped,below] {$s$=1}},
                    edge node={node [yshift = -0.23cm] {$>$}}](q_2);
\draw[->] (q_1) to [looseness=1.4, bend left=45, 
                    edge node={node [sloped,below] {$s$=2}},
                    edge node={node [yshift = -0.23cm] {$>$}}](q_3);
\draw[->] (q_1) to [looseness=1.4, bend left=45, 
                    edge node={node [sloped,below] {$s$=3}},
                    edge node={node [yshift = -0.23cm] {$>$}}](q_4);
\draw[->] (q_2) to [looseness=1.4, bend left=45, 
                    edge node={node [sloped,below] {$s$=1}},
                    edge node={node [yshift = -0.23cm] {$>$}}](q_3);
\draw[->] (q_2) to [looseness=1.4, bend left=45, 
                    edge node={node [sloped,below] {$s$=2}},
                    edge node={node [yshift = -0.23cm] {$>$}}](q_4);
\draw[->] (q_2) to [looseness=1.4, bend left=45, 
                    edge node={node [sloped,below] {$s$=3}},
                    edge node={node [yshift = -0.23cm] {$>$}}](q_5);
\draw[->] (q_3) to [looseness=1.4, bend left=45, 
                    edge node={node [sloped,below] {$s$=1}},
                    edge node={node [yshift = -0.23cm] {$>$}}](q_4);
\draw[->] (q_3) to [looseness=1.4, bend left=45, 
                    edge node={node [sloped,below] {$s$=2}},
                    edge node={node [yshift = -0.23cm] {$>$}}](q_5);
\draw[->] (q_4) to [looseness=1.4, bend left=45, 
                    edge node={node [sloped,below] {$s$=1}},
                    edge node={node [yshift = -0.23cm] {$>$}}](q_5);
\draw[->] (q_1) to [bend right=40, 
                    edge node={node [sloped,above] {$loss$}},
                    edge node={node [yshift = -0.23cm] {$>$}}](q_5);
\draw[->] (q_2) to [bend right=40, 
                    edge node={node [sloped,above] {$loss$}},
                    edge node={node [yshift = -0.23cm] {$>$}}](q_5);
\draw[->] (q_3) to [bend right=40, 
                    edge node={node [sloped,above] {$loss$}},
                    edge node={node [yshift = -0.23cm] {$>$}}](q_5);
\draw[->] (q_4) to [bend right=40, 
                    edge node={node [sloped,above] {$loss$}},
                    edge node={node [yshift = -0.23cm] {$>$}}](q_5);
\draw[dashed, <-] (q_0) |- +(0,-2.8) -| 
                  node[pos=0.25] {$feedback$} 
                  node[pos=0.25, yshift = -0.23cm] {$<$} (q_5.south);
\end{tikzpicture}
\end{figure}

Figure \ref{fig_flow1}
shows the graph $G$  for a machine with capacity $B=5$, and a set of jobs $\{1,2,3,4,5\}$ with respective sizes $\{s_1=3,s_2=2,s_3=2,s_4=1,s_5=1\}$. 


In our modeling approach, each unit flow in graph $G$, going from node $0$ to node $B$,  represents  the  configuration of a batch in the solution of the problem. Moreover, a unit flow from node $B$ to node $ 0 $, going through the \emph{feedback arc},  represents a solution using  a single batch, and flows with more than one unit on the \emph{feedback arc}, represent the use of several batches. 

As all flows are non-negative and $G$ is acyclic once the \emph{feedback arc} is excluded, it is possible to decompose the multiple flow on the \emph{feedback arc} into several unit flows, each corresponding to an oriented path from node  $0$ to node $B$. Each unit flow represents the  configuration of a different batch and the flow on the \emph{feedback arc} indicates the number of batches used in the solution. 

\begin{figure}[!ht]
\begin{center}
\subfigure[]{
    \label{fig_flowDec_a}
    \begin{tikzpicture}[shorten >=1pt,node distance=2.5cm,auto]
\node[circle, draw] (q_0) {$0$};
\node[circle, draw] (q_1) [right of=q_0] {$1$};
\node[circle, draw] (q_2) [right of=q_1] {$2$};
\node[circle, draw] (q_3) [right of=q_2] {$3$};
\node[circle, draw] (q_4) [right of=q_3] {$4$};
\node[circle, draw] (q_5) [right of=q_4] {$5$};
\draw[->] (q_0) to [looseness=1.4, bend left=45, 
                    edge node={node [sloped,below] {1}},
                    edge node={node [yshift = -0.23cm] {$>$}}](q_2);
\draw[->] (q_0) to [looseness=1.4, bend left=45, 
                    edge node={node [sloped,below] {1}},
                    edge node={node [yshift = -0.23cm] {$>$}}](q_3);
\draw[->] (q_2) to [looseness=1.4, bend left=45, 
                    edge node={node [sloped,below] {1}},
                    edge node={node [yshift = -0.23cm] {$>$}}](q_4);
\draw[->] (q_3) to [looseness=1.4, bend left=45, 
                    edge node={node [sloped,below] {1}},
                    edge node={node [yshift = -0.23cm] {$>$}}](q_4);
\draw[->] (q_4) to [looseness=1.4, bend left=45, 
                    edge node={node [sloped,below] {1}},
                    edge node={node [yshift = -0.23cm] {$>$}}](q_5);
\draw[->] (q_4) to [bend right=40, 
                    edge node={node [sloped,above] {1}},
                    edge node={node [yshift = -0.23cm] {$>$}}](q_5);
\draw[dashed, <-] (q_0) |- +(0,-1.5) -| 
                  node[pos=0.25] {2} 
                  node[pos=0.25, yshift = -0.23cm] {$<$} (q_5.south);
\end{tikzpicture}
}
\subfigure[]{
    \label{fig_flowDec_b}
    \begin{tikzpicture}[shorten >=1pt,node distance=2.5cm,auto]
\node[circle, draw] (q_0) {$0$};
\node[circle, draw] (q_1) [right of=q_0] {$1$};
\node[circle, draw] (q_2) [right of=q_1] {$2$};
\node[circle, draw] (q_3) [right of=q_2] {$3$};
\node[circle, draw] (q_4) [right of=q_3] {$4$};
\node[circle, draw] (q_5) [right of=q_4] {$5$};
\draw[->] (q_0) to [looseness=1.4, bend left=45, 
                    edge node={node [sloped,below] {$s_1=3$}},
                    edge node={node [yshift = -0.23cm] {$>$}}](q_3);
\draw[->] (q_3) to [looseness=1.4, bend left=45, 
                    edge node={node [sloped,below,yshift = -0.2cm] {$s_4=1$}},
                    edge node={node [yshift = -0.23cm] {$>$}}](q_4);
\draw[->] (q_4) to [looseness=1.4, bend left=45, 
                    edge node={node [sloped,below,yshift = -0.2cm] {$s_5=1$}},
                    edge node={node [yshift = -0.23cm] {$>$}}](q_5);
\draw[dashed, <-] (q_0) |- +(0,-1.5) -| 
                  node[pos=0.25] {$feedback$} 
                  node[pos=0.25, yshift = -0.23cm] {$<$} (q_5.south);
\end{tikzpicture}
}
\subfigure[]{
    \label{fig_flowDec_c}
    \begin{tikzpicture}[shorten >=1pt,node distance=2.5cm,auto]
\node[circle, draw] (q_0) {$0$};
\node[circle, draw] (q_1) [right of=q_0] {$1$};
\node[circle, draw] (q_2) [right of=q_1] {$2$};
\node[circle, draw] (q_3) [right of=q_2] {$3$};
\node[circle, draw] (q_4) [right of=q_3] {$4$};
\node[circle, draw] (q_5) [right of=q_4] {$5$};
\draw[->] (q_0) to [looseness=1.4, bend left=45, 
                    edge node={node [sloped,below] {$s_2=2$}},
                    edge node={node [yshift = -0.23cm] {$>$}}](q_2);
\draw[->] (q_2) to [looseness=1.4, bend left=45, 
                    edge node={node [sloped,below] {$s_3=2$}},
                    edge node={node [yshift = -0.23cm] {$>$}}](q_4);
\draw[->] (q_4) to [bend right=40, 
                    edge node={node [sloped,above] {$loss$}},
                    edge node={node [yshift = -0.23cm] {$>$}}](q_5);
\draw[dashed, <-] (q_0) |- +(0,-1.5) -| 
                  node[pos=0.25] {$feedback$} 
                  node[pos=0.25, yshift = -0.23cm] {$<$} (q_5.south);
\end{tikzpicture}
}
\caption{(a) Feasible solution for problem represented in Figure \ref{fig_flow1}, with multi-flow decomposed into (b) and (c).} 
\label{fig_flowDec}
\end{center}
\end{figure}
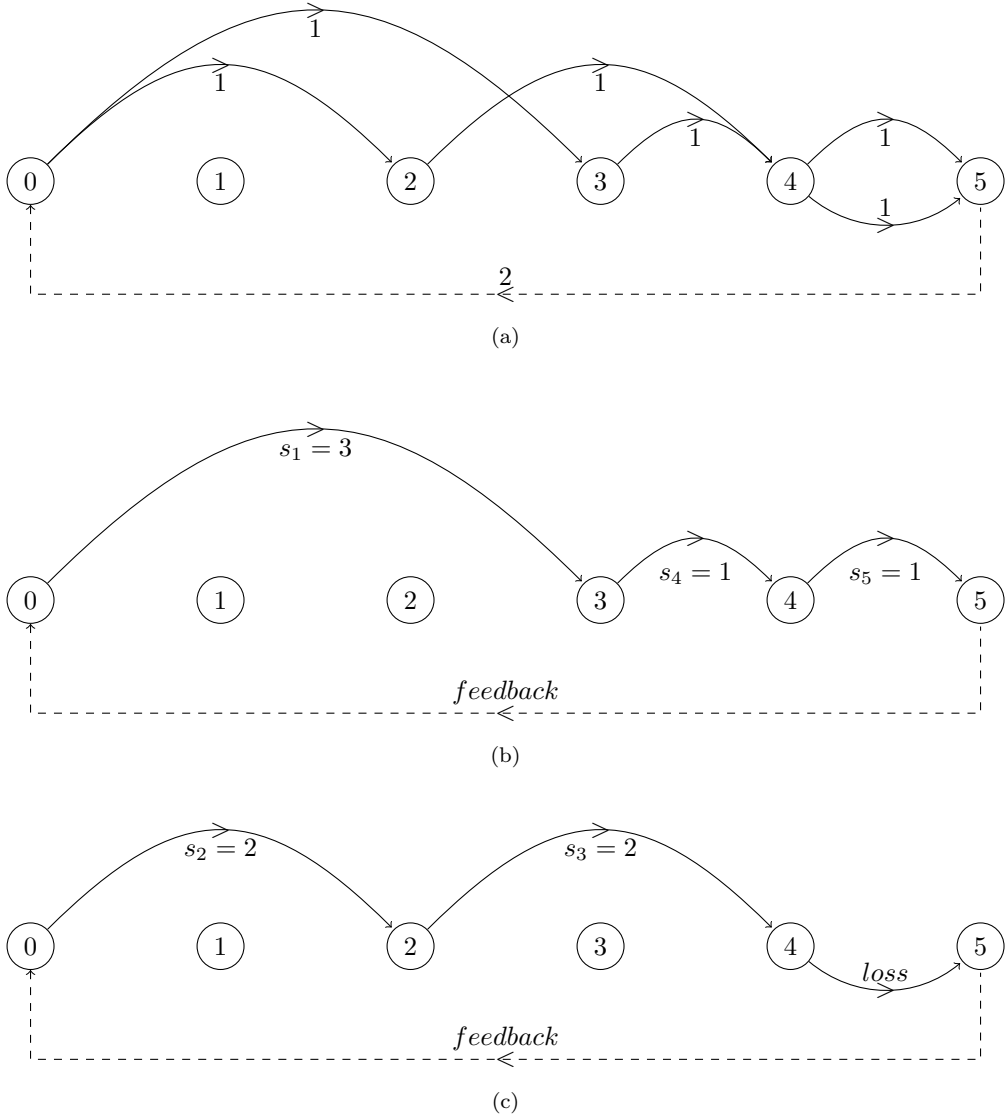

Figure \ref{fig_flowDec_a} depicts a feasible solution for the problem represented by the graph $G$ in Figure \ref{fig_flow1}. The multiple-flow graph of two units is decomposed into two other graphs in Figures \ref{fig_flowDec_b} and \ref{fig_flowDec_c}, each with a unit flow. They represent  two batches  $ 1 $ and $ 2 $, such that jobs $\{1, 4, 5 \} $ are assigned to batch $1$  (Figure \ref{fig_flowDec_b}) and jobs $\{2, 3\}$ are assigned to batch $2$ (Figure \ref{fig_flowDec_c}). Arcs 
with  null value  are not shown.

The graph $G$ defined above, is then replicated for each distinct processing time of the problem in our modeling approach. Each replicated graph will be referred to, as an arc-flow structure for our problem. We define $P:=\{P_1,\ldots,P_{\delta}\}$ as the set with all the distinct processing times among all jobs, and $T:=\{1,\ldots,\delta\}$ as the set of indexes corresponding to the arc-flow structures in the problem formulation. We consider $P_1< P_2<  \ldots < P_\delta$. 
The arc-flow structure $t\in T$, will have the corresponding processing time fixed at $P_t$, meaning that all batches represented in the arc-flow structure $t$ will be considered in our formulation to have processing time $P_t$.  Only jobs $j$ with processing time $p_j \leq P_t $ will be allowed to be assignment to this arc-flow structure. 
Moreover, the minimization of the makespan in the objective function will assure that all batches represented in the arc-flow structure $t$ in the optimal solution of the problem, have at least one job with processing time $P_t$ assigned to it. The number of batches corresponding to each arc-flow structure is indicated by the value of the flow in the \emph{feedback arc}. 
As each arc-flow structure can represent as many batches as needed, it is possible to represent all solutions of problem  $1|s_j,B|C_{\max}$ with our arc-flow structures. Finally, we note that, in the optimal solution of the problem, the maximum number of batches represented in the arc-flow structure $t$ is given by the number of jobs with processing time $P_t$. This bound will only be attained when each job with processing time $P_t$ is assigned to a different batch.

\begin{figure}[!ht]
\begin{center}
\def\scale{0.35}
\newcommand{\GenerateXAxis}[0]{
    \draw [->, thick] (0,0-0.5) node [yshift = -0.5cm] {0} -- node [yshift = -0.8cm] {\NameXAxis} (\TimeMaxValue+0.5,0-0.5) node [yshift = -0.5cm] {$C_{max}$};
    \foreach \i in {1,...,\TimeMaxValue}
    {
        \draw[-] (\i,-0.5) -- (\i,-1);
    }
    \draw[-, thick] (0,-0.5) -- (0,-1.18);
    \draw[-, thick] (\TimeMaxValue,-0.5) -- (\TimeMaxValue,-1.18);
}

\newcommand{\GenerateSizeAxe}[1]{
    \draw [->, thick] (0-0.5,\SizeBeginValueY + 0) node [xshift = -0.5cm] {0} -- node [xshift = -0.2cm] (#1) {Size} (0-0.5,\SizeBeginValueY + \SizeMaxValue+0.5) node [xshift = -0.5cm, yshift = -0.2cm] {$B$};
    
    \foreach \i in {1,...,\SizeMaxValue}
    {
        \draw[-] (-0.5,\SizeBeginValueY + \i) -- (-1,\SizeBeginValueY + \i);
    }
    \draw[-, thick] (-0.5,\SizeBeginValueY) -- (-1.18,\SizeBeginValueY);
    \draw[-, thick] (-0.5,\SizeBeginValueY+\SizeMaxValue) -- (-1.18,\SizeBeginValueY+\SizeMaxValue);
}

\newcommand{\CreateNewJob}[4]{
    \node[fit={#4 ($#4 +
    #1
    $)}, inner sep=0.3pt, draw=black, shade, left color=mycolor, right color=white, label=center:#3] (#2) {};
}

\newcommand{\CreateNewJobDashed}[4]{
    \node[fit={#4 ($#4 +
    #1
    $)}, inner sep=0.3pt, draw=black, dashed, shade, left color=white, right color=white, label=center:#3] (#2) {};
}

\begin{tikzpicture}[scale=\scale, shorten >=1pt,node distance=2.5cm,auto]
\definecolor{mycolor}{RGB}{220,220,220}
\tikzmath{
    \SizeMaxValue = 5;
    \TimeMaxValue = 24;
    \SizeBeginValueY = 0;
    \DistanceBtweenGraphs = 0;
    \BatchRelease = 0;
    \TopOfGraph = \SizeMaxValue+1;
}

\draw [->, thick] node [yshift = -0.8cm] {0} (0,0-0.5) -- node [yshift = -0.8cm] {Time} (\TimeMaxValue+0.5,0-0.5) node [yshift = -0.6cm] {$C_{max}$};
\foreach \i in {1,...,\TimeMaxValue}
{
    \draw[-] (\i,-0.5) -- (\i,-1);
}
\draw[-, thick] (0,-0.5) -- (0,-1.18);
\draw[-, thick] (\TimeMaxValue,-0.5) -- (\TimeMaxValue,-1.18);

\GenerateSizeAxe{batch1}

\draw[-, thick, dashed] (\TimeMaxValue,-0.5) -- (\TimeMaxValue,\TopOfGraph);

\draw[-, thick, dashed, , opacity=0.4] (0,-0.5) -- (0,\TopOfGraph);
\draw[-, thick, dashed, , opacity=0.4] (3,-0.5) -- (3,\TopOfGraph);
\draw[-, thick, dashed, , opacity=0.4] (6,-0.5) -- (6,\TopOfGraph);
\draw[-, thick, dashed, , opacity=0.4] (10,-0.5) -- (10,\TopOfGraph);
\draw[-, thick, dashed, , opacity=0.4] (14,-0.5) -- (14,\TopOfGraph);
\draw[-, thick, dashed, , opacity=0.4] (19,-0.5) -- (19,\TopOfGraph);




\tikzmath{
    \SizeBeginValueY = 0;
    \ReleaseTime = 0;
}
\CreateNewJob{(3,2)}{job_3}{3}{(\ReleaseTime,\SizeBeginValueY)};
\CreateNewJob{(3,2)}{job_1}{1}{($(job_3.north west) + (1pt,0)$)};
\draw[<->] ($(job_3.south west) + (0,\TopOfGraph)$) -- node [above] {$P_1=3$} ($(job_3.south east) + (0,\TopOfGraph)$);

\CreateNewJob{(3,2)}{job_5}{5}{($(job_3.south east)+(-1pt,1pt)$)};
\CreateNewJob{(3,1)}{job_2}{2}{($(job_5.north west) + (1pt,0)$)};
\CreateNewJob{(3,1)}{job_4}{4}{($(job_2.north west) + (1pt,0)$)};
\draw[<->] ($(job_5.south west) + (0,\TopOfGraph)$) -- node [above] {$P_1=3$} ($(job_5.south east) + (0,\TopOfGraph)$);

\CreateNewJob{(4,1)}{job_9}{9}{($(job_5.south east)+(-1pt,1pt)$)};
\CreateNewJob{(4,3)}{job_6}{6}{($(job_9.north west) + (1pt,0)$)};
\CreateNewJob{(4,1)}{job_8}{8}{($(job_6.north west) + (1pt,0)$)};
\draw[<->] ($(job_9.south west) + (0,\TopOfGraph)$) -- node [above] {$P_2=4$} ($(job_9.south east) + (0,\TopOfGraph)$);

\CreateNewJob{(4,3)}{job_11}{11}{($(job_9.south east)+(-1pt,1pt)$)};
\CreateNewJob{(4,2)}{job_10}{10}{($(job_11.north west) + (1pt,0)$)};
\draw[<->] ($(job_11.south west) + (0,\TopOfGraph)$) -- node [above] {$P_2=4$} ($(job_11.south east) + (0,\TopOfGraph)$);

\CreateNewJob{(5,2)}{job_14}{14}{($(job_11.south east)+(-1pt,1pt)$)};
\CreateNewJob{(5,3)}{job_12}{12}{($(job_14.north west) + (1pt,0)$)};
\draw[<->] ($(job_14.south west) + (0,\TopOfGraph)$) -- node [above] {$P_3=5$} ($(job_14.south east) + (0,\TopOfGraph)$);

\CreateNewJob{(5,2)}{job_15}{15}{($(job_14.south east)+(-1pt,1pt)$)};
\CreateNewJob{(4,2)}{job_7}{7}{($(job_15.north west) + (1pt,0)$)};
\CreateNewJob{(5,1)}{job_13}{13}{($(job_7.north west) + (1pt,0)$)};
\draw[<->] ($(job_15.south west) + (0,\TopOfGraph)$) -- node [above] {$P_3=5$} ($(job_15.south east) + (0,\TopOfGraph)$);

\end{tikzpicture}


\caption{Solution with $C_{max}=24$ and six batches.} 
\label{fig_NewExampleSingle}
\end{center}
\end{figure}
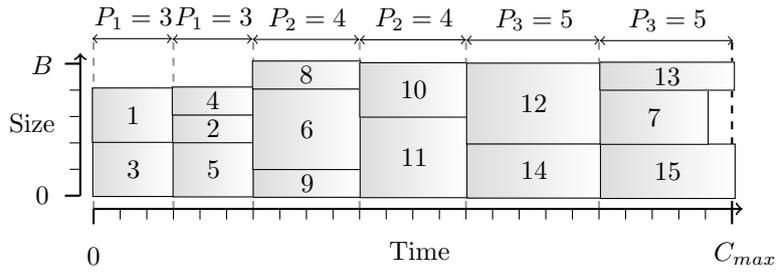

\begin{figure}[!ht]
\begin{center}
\subfigure[
Representation of batches  with processing time $P_1=3$
]{
    \label{fig_NewExampleSingle_a}
    \begin{tikzpicture}[shorten >=1pt,node distance=2.5cm,auto]
\node[circle, draw] (q_0) {$0$};
\node[circle, draw] (q_1) [right of=q_0] {$1$};
\node[circle, draw] (q_2) [right of=q_1] {$2$};
\node[circle, draw] (q_3) [right of=q_2] {$3$};
\node[circle, draw] (q_4) [right of=q_3] {$4$};
\node[circle, draw] (q_5) [right of=q_4] {$5$};
\draw[->] (q_0) to [looseness=1.4, bend left=45, 
                    edge node={node [sloped,below] {2}},
                    edge node={node [yshift = -0.23cm] {$>$}}](q_2);
\draw[->] (q_2) to [looseness=1.4, bend left=45, 
                    edge node={node [sloped,below] {1}},
                    edge node={node [yshift = -0.23cm] {$>$}}](q_3);
\draw[->] (q_2) to [looseness=1.4, bend left=45, 
                    edge node={node [sloped,below] {1}},
                    edge node={node [yshift = -0.23cm] {$>$}}](q_4);
\draw[->] (q_3) to [looseness=1.4, bend left=45, 
                    edge node={node [sloped,below] {1}},
                    edge node={node [yshift = -0.23cm] {$>$}}](q_4);
\draw[->] (q_4) to [bend right=40, 
                    edge node={node [sloped,above] {2}},
                    edge node={node [yshift = -0.23cm] {$>$}}](q_5);
\draw[dashed, <-] (q_0) |- +(0,-1.5) -| 
                  node[pos=0.25] {2} 
                  node[pos=0.25, yshift = -0.23cm] {$<$} (q_5.south);
\end{tikzpicture}
}
\subfigure[
Representing of batches  with processing time $P_2=4$
]{
    \label{fig_NewExampleSingle_b}
    \begin{tikzpicture}[shorten >=1pt,node distance=2.5cm,auto]
\node[circle, draw] (q_0) {$0$};
\node[circle, draw] (q_1) [right of=q_0] {$1$};
\node[circle, draw] (q_2) [right of=q_1] {$2$};
\node[circle, draw] (q_3) [right of=q_2] {$3$};
\node[circle, draw] (q_4) [right of=q_3] {$4$};
\node[circle, draw] (q_5) [right of=q_4] {$5$};
\draw[->] (q_0) to [looseness=1.4, bend left=45, 
                    edge node={node [sloped,below] {1}},
                    edge node={node [yshift = -0.23cm] {$>$}}](q_1);
\draw[->] (q_0) to [looseness=1.4, bend left=45, 
                    edge node={node [sloped,below] {1}},
                    edge node={node [yshift = -0.23cm] {$>$}}](q_3);
\draw[->] (q_1) to [looseness=1.4, bend left=45, 
                    edge node={node [sloped,below] {1}},
                    edge node={node [yshift = -0.23cm] {$>$}}](q_4);
\draw[->] (q_3) to [looseness=1.4, bend left=45, 
                    edge node={node [sloped,below] {1}},
                    edge node={node [yshift = -0.23cm] {$>$}}](q_5);
\draw[->] (q_4) to [looseness=1.4, bend left=45, 
                    edge node={node [sloped,below] {1}},
                    edge node={node [yshift = -0.23cm] {$>$}}](q_5);
\draw[dashed, <-] (q_0) |- +(0,-1.5) -| 
                  node[pos=0.25] {2} 
                  node[pos=0.25, yshift = -0.23cm] {$<$} (q_5.south);
\end{tikzpicture}
}
\subfigure[
Representation of batches  with processing time $P_3=5$
]{
    \label{fig_NewExampleSingle_c}
    \begin{tikzpicture}[shorten >=1pt,node distance=2.5cm,auto]
\node[circle, draw] (q_0) {$0$};
\node[circle, draw] (q_1) [right of=q_0] {$1$};
\node[circle, draw] (q_2) [right of=q_1] {$2$};
\node[circle, draw] (q_3) [right of=q_2] {$3$};
\node[circle, draw] (q_4) [right of=q_3] {$4$};
\node[circle, draw] (q_5) [right of=q_4] {$5$};
\draw[->] (q_0) to [looseness=1.4, bend left=45, 
                    edge node={node [sloped,below] {2}},
                    edge node={node [yshift = -0.23cm] {$>$}}](q_2);
\draw[->] (q_2) to [looseness=1.4, bend left=45, 
                    edge node={node [sloped,below] {1}},
                    edge node={node [yshift = -0.23cm] {$>$}}](q_4);
\draw[->] (q_2) to [looseness=1.4, bend left=45, 
                    edge node={node [sloped,below] {1}},
                    edge node={node [yshift = -0.23cm] {$>$}}](q_5);
\draw[->] (q_4) to [looseness=1.4, bend left=45, 
                    edge node={node [sloped,below] {1}},
                    edge node={node [yshift = -0.23cm] {$>$}}](q_5);
\draw[dashed, <-] (q_0) |- +(0,-1.5) -| 
                  node[pos=0.25] {2} 
                  node[pos=0.25, yshift = -0.23cm] {$<$} (q_5.south);
\end{tikzpicture}
}
\caption{Arc-flow structures representing the solution in Figure \ref{fig_NewExampleSingle}.} 
\label{fig_NewExampleSingleFlow}
\end{center}
\end{figure}
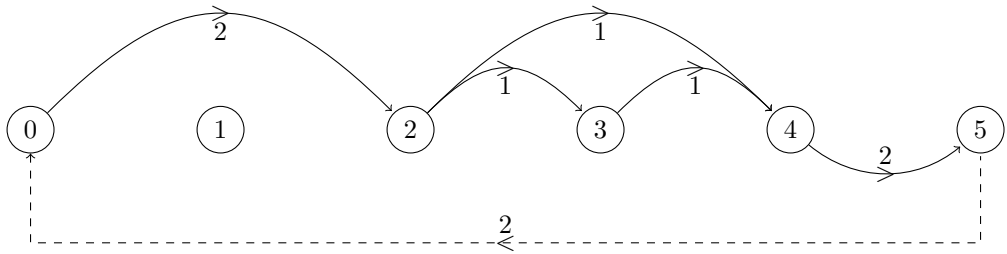
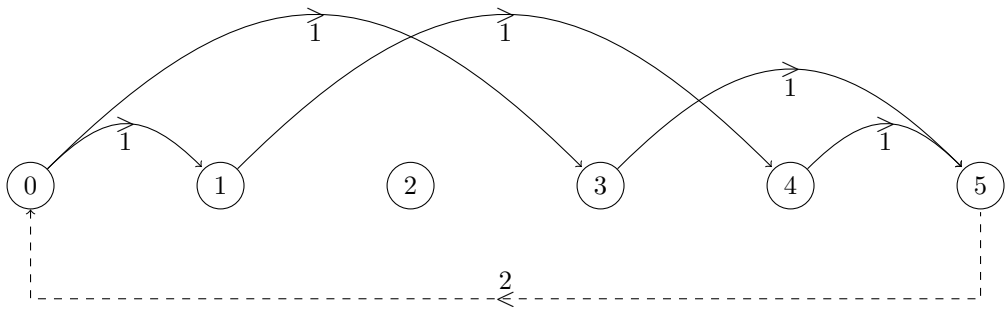
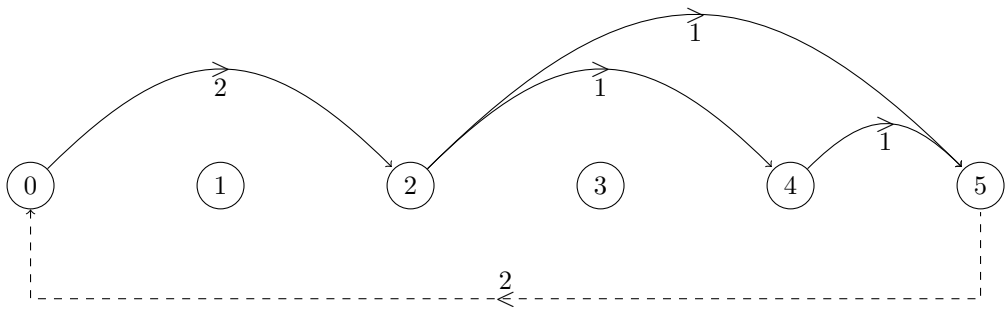

In Figure \ref{fig_NewExampleSingle}, we illustrate a feasible solution for an instance of problem  $1|s_j,B|C_{\max}$ with fifteen jobs. The processing times for the jobs are either  3, 4, or 5. They are grouped into 
six batches, two of them per processing time. 
We want to emphasize here, that the number of  the arc-flow structures used in our formulation is given by the number of distinct processing times on the instance. Therefore, for this example,  only three structures are used to represent the feasible solutions. Particularly, the solution in Figure  \ref{fig_NewExampleSingle} is represented by the three structures  in Figure \ref{fig_NewExampleSingleFlow},  one for each distinct processing time. 

The graph $G$ and the arc-flow structures introduced  above form together the core of the arc-flow formulation for problem  $1|s_j,B|C_{\max}$. It is important then to analyze what are the parameters of the problem that interfere in the size and density of $G$ and also in the number of arc-flow structures that are needed to represent the solution.  These are the parameters that will interfere in the size of our problem formulation. We emphasize them in the following remark.

\begin{remark}
\label{sizeform}
Consider problem  $1|s_j,B|C_{\max}$, with machine capacity $B\in\mathbb{Z}$. There are a set of jobs to be assigned to batches, with $\delta$ distinct processing times $\{P_1,P_2,\ldots,P_\delta\}$ and  $\theta$ distinct sizes $\{S_1,S_2,\ldots,S_\theta\}$.
Consider the graph $G=(V,A)$ and the arc-flow structures defined as above to represent a feasible solution to  $1|s_j,B|C_{\max}$. We have $|V|=B+1$,  $|A|=\theta+(\theta+1) B-\sum_{i=1}^{\theta}S_i$, and the number of arc-flow structures needed to represent all feasible solutions equal to $\delta$. 
\end{remark}

We observe from Remark \ref{sizeform} that  the size of of the arc-flow formulation that we will present for problem $1|s_j,B|C_{\max}$, increases with the number of distinct sizes and processing times among the jobs, but  does not increase with the number of jobs. This makes the formulation particularly suitable for large instances where there are  multiple jobs with similar characteristics, specifically with  equal  size and processing time.

 
\subsection{Mapping arc flows into job assignments and breaking symmetries}
\label{subsecmap}
 
The solution represented by the arc-flow structure described above, does not specify exactly which are the jobs assigned to each batch, but only the number of jobs with each size that are assigned to it.  
For example, if we decompose the two-units  multi-flow structure in Figure \ref{fig_NewExampleSingle_a} into two unit-flow structures, we can easily see that it represents two batches with processing time $P_1=3$. One of them has two jobs of size $s=2$ and the other has  two jobs of size $s=1$ and one job of size $s=2$. As shown in Figure \ref{fig_NewExampleSingle}, the instance considered has  three jobs of size $s=2$ and processing time $P_1=3$, jobs  1,3, and 5. Therefore, each one of them could be assigned to the latter batch, generating three equivalent solutions. 

Moreover, the arc-flow solution determines the number of batches with each processing time, but defines no processing ordering for them. As in  problem  $1|s_j,B|C_{\max}$, release times are not considered for the jobs, any permutation in the ordering of the batches would also lead to equivalent solutions.

Frequently, mathematical programming formulations for combinatorial optimization problems, make distinction between equivalent solutions by indexing them in different ways. These are called symmetric solutions, and their existence in the feasible set of the problems tends to worsen the performance of branch-and-bound algorithms for wasting effort evaluating equivalent solution unnecessarily. In \cite{Trindade2018}, the authors point out that this is the case of formulation (MILP$_1$), and show that the symmetry generated by exchanging the indexes assigned to the batches in the solution of (MILP$_1$)  is broken in (MILP$_1^+$). We note, however, that the symmetry generated by exchanging the indexes assigned to jobs with the same sizes and processing times is still not broken in (MILP$_1^+$). 

The arc-flow formulation have neither of these  symmetries in its feasible set, and with a simple mapping, we can determine a feasible assignment of jobs to batches, according to the arc-flow optimal solution.
The mapping of arcs into jobs  is easily performed in polynomial  time, since each \emph{job arc} $(i,j)$ in the arc-flow structure corresponding to processing time $P_t$, can be mapped into any job of size $j-i$ and processing time not greater than $P_t$. 
In order to guarantee that a feasible solution is obtained, we should only map the arcs into jobs starting from the structures with the smallest processing time and continue in increasing order of processing time.  Concerning structures with the same processing times, the ordering for the mapping among them is not relevant.

Finally, we should observe that, although we have eliminated symmetries contained in the feasible sets of (MILP$_1$) and  (MILP$_1^+$) with the arc-flow formulation, it still contains symmetric solutions in its feasible set, related to the ``positions'' occupied by the jobs in the batches, as determined by the definition of the variables of the formulation presented in the following.

\subsection{Arc reduction}

Some rules can be set to decrease the number of arcs in graph $G$ and the number of variables in our arc-flow formulation.
The set  $A^L $  of  \emph{loss arcs} defined in (\ref{arc_L}), represents empty spaces that occur at the end of the batch.
We first note that our modeling does not represent empty spaces between jobs or at the beginning of the batch, which avoids symmetric solutions on the feasible set of the problem. Because of that, two rules can be created to eliminate arcs that will never be used in the solution of the problem: Rule \ref{rulearc1} eliminates \emph{job arcs} that cannot be used in a valid flow starting at node 0, and Rule \ref{rulearc2} eliminates  \emph{loss arcs}  incident to nodes that are not incident to any \emph{job arc}.

\begin{rules}
Only \emph{job arcs} $(i, j) \in A^J$  that belong to at least one continuous flow  starting at node 0 can have a positive flow in the solution of the problem, and therefore, all the others may be eliminated from graph $G$. The \emph{job arcs} that belong to at least one continuous flow can be selected by the following steps:\\
\begin{enumerate}
\vspace{-0.22in}
\item arc $(0,j)$ is selected, for all $j$, such that $(0,j)\in A^J$;\\
\vspace{-0.25in}
\item arc $(i, j)\in A^J$ is selected, if an arc $(k,i)$ has been previously selected for some node $k$;\\
\vspace{-0.25in}
\item repeat step (2) until no arc is selected.\\
\vspace{-0.25in}
\end{enumerate}
\label{rulearc1}
\end{rules} 

\begin{rules}
Only \emph{loss arcs }  $(i,B) \in A^L $  that are incident to a node $i$, which is incident to a remaining \emph{job arc} after the application of Rule \ref{rulearc1}, can have a positive flow in the solution of the problem, and therefore, all the others may be eliminated from graph $G$.
\label{rulearc2}
\end{rules} 

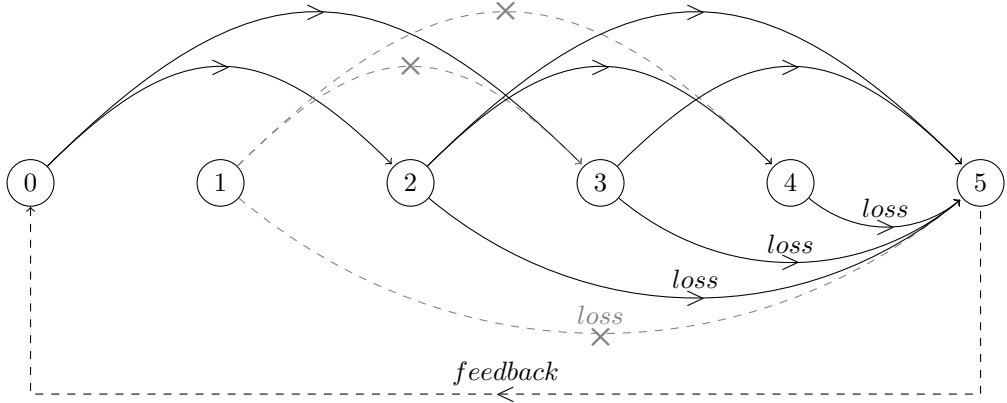
\begin{figure}[!ht]
	\centering
	\begin{tikzpicture}[shorten >=1pt,node distance=2.5cm,auto]
\node[circle, draw] (q_0) {$0$};
\node[circle, draw] (q_1) [right of=q_0] {$1$};
\node[circle, draw] (q_2) [right of=q_1] {$2$};
\node[circle, draw] (q_3) [right of=q_2] {$3$};
\node[circle, draw] (q_4) [right of=q_3] {$4$};
\node[circle, draw] (q_5) [right of=q_4] {$5$};
\draw[->] (q_0) to [looseness=1.4, bend left=45, 
                    edge node={node [sloped,below] {}},
                    edge node={node [yshift = -0.23cm] {$>$}}](q_2);
\draw[->] (q_0) to [looseness=1.4, bend left=45, 
                    edge node={node [sloped,below] {}},
                    edge node={node [yshift = -0.23cm] {$>$}}](q_3);
\draw[->, gray, dashed] (q_1) to [looseness=1.4, bend left=45, 
                    edge node={node [sloped,below] {}},
                    edge node={node [sloped,above] {}},
                    edge node={node [yshift = -0.23cm] {\ding{53}}}
                    ](q_3);
\draw[->, gray, dashed] (q_1) to [looseness=1.4, bend left=45, 
                    edge node={node [sloped,below] {}},
                    edge node={node [sloped,above] {}},
                    edge node={node [yshift = -0.23cm] {\ding{53}}}
                    ](q_4);
\draw[->] (q_2) to [looseness=1.4, bend left=45, 
                    edge node={node [sloped,below] {}},
                    edge node={node [yshift = -0.23cm] {$>$}}](q_4);
\draw[->] (q_2) to [looseness=1.4, bend left=45, 
                    edge node={node [sloped,below] {}},
                    edge node={node [yshift = -0.23cm] {$>$}}](q_5);
\draw[->] (q_3) to [looseness=1.4, bend left=45, 
                    edge node={node [sloped,below] {}},
                    edge node={node [yshift = -0.23cm] {$>$}}](q_5);
%
\draw[->, gray, dashed] (q_1) to [bend right=40, 
                    edge node={node [sloped,above] {$loss$}},
                    edge node={node [yshift = -0.29cm] {\ding{53}}}
                    ](q_5);
\draw[->] (q_2) to [bend right=40, 
                    edge node={node [sloped,above] {$loss$}},
                    edge node={node [yshift = -0.23cm] {$>$}}](q_5);
\draw[->] (q_3) to [bend right=40, 
                    edge node={node [sloped,above] {$loss$}},
                    edge node={node [yshift = -0.23cm] {$>$}}](q_5);
\draw[->] (q_4) to [bend right=40, 
                    edge node={node [sloped,above] {$loss$}},
                    edge node={node [yshift = -0.23cm] {$>$}}](q_5);
\draw[dashed, <-] (q_0) |- ($ (q_0) - (0,2.8) $) -| 
                  node[pos=0.25] {$feedback$} 
                  node[pos=0.25, yshift = -0.23cm] {$<$} (q_5.south);
\end{tikzpicture}
	\caption{Example of arc reduction. }
	\label{fig_flowArcReduction}
\end{figure}

{\color{red} 

}
Figure \ref{fig_flowArcReduction} illustrates the application of Rules \ref{rulearc1}--\ref{rulearc2}, where 2 \emph{job arcs} and  1 \emph{loss arc} are removed from graph $G$.

\subsection{Defining bounds for the  variables in the formulation}

In   the following, we define bounds on  the values of the flows in each arc of the arc-flow structures, and therefore bounds for the variables in the  model poposed . 

\begin{rules}
The flow on each arc $(i,j)\in A$ that belongs to the arc-flow structure corresponding to processing time $p$ cannot be greater than the number of jobs with processing time $p$. This upper bound is more specifically given by
\begin{align}\nonumber
\sum_{\substack{ k \in J: p_k=p}} 1.
\end{align}
\label{rulearc3}
\end{rules} 

\begin{rules}
The flow on each  \emph{job arc} $(i,j)\in A^J$ that belongs to the arc-flow structure corresponding to processing time $p$ cannot be greater than the number of jobs of size $j-i$ and processing time not greater than $p$. This upper bound is more specifically given by
    \begin{align}\nonumber
	\sum_{\substack{ k \in J: s_k=j-i,\\ p_k \leq p}} 1.
	\end{align}
\label{rulearc4}
\end{rules}

  \subsection{Arc-flow formulation}\label{sec:form}

In this section, we formulate problem $1|s_j,B|C_{\max}$ as the problem of determining the flows from node 0 to node $B$ on the arc-flow structures associated to the distinct processing times of the jobs. The objective is to minimize  the  sum of the flows on the structures, weighted by their associated processing times. We emphasize that the flow on the structure associated to the processing time $P_t$ indicates the number of batches with processing time $P_t$ in the solution, therefore, the objective function of our model represents the makespan of the batch processing machine.

Besides the set of variables and constraints that determine the flows on the arcs of the structures and ensure the flow conservation on the nodes, a set of variables and constraints is also required to ensure that all jobs are assigned to an arc-flow structure (or equivalently, to a batch), and that each  structure considers the correct number of jobs of each size that are available to be assigned to it. 

We recall that our approach requires the arc-flow structures to be sorted by increasing processing time, i.e., the structure $t$ is associate to the processing time $P_t$, and $P_1<P_2<\ldots<P_\delta$, where $\delta$  is the number of distinct processing times among all jobs. Let us then suppose that the arc-flow structure $t$  has $\beta$ jobs of size $S_\ell$ still available to be assigned to it. These are all the jobs of size $S_\ell$ and processing time not greater than $P_t$, that have not already been assigned to an arc-flow structure with corresponding processing time smaller than $P_t$. A variable $z_{\ell,t}$ is defined in our model to represent how many of these $\beta$ jobs of size $S_\ell$  were not assigned to  the arc-flow structure  $t$. 
As $\beta$ is the total number of jobs available and assuming that $\alpha$ is the number of jobs that were assigned to structure $t$, the variable $z_{\ell,t}$  assumes value $\beta-\alpha$. This variable will  ``offer'' the remaining  $\beta-\alpha$ jobs of size $S_\ell$ to the next arc-flow  structure, 
$t+1$, corresponding to next processing time of the problem, greater than $P_t$. In the end, all jobs should be assigned to a structure, and, therefore, we should have  $z_{\ell,\delta}=0$.

Our arc-flow formulation is presented next.

\vspace{0.1in}

	\noindent
	\textbf{Sets and parameters}\\
	$B$ : machine capacity.\\
	$\delta$ :  number of distinct processing times among all jobs.\\
	$P:=\{P_1,\ldots, P_{\delta}\}$ :  set of distinct processing times of jobs, where $P_1<P_2<\ldots<P_\delta$.\\
	$T:=\{1,\ldots,\delta\}$ : set of indexes of the arc-flow structures.\\
	$\theta$ :  number of distinct sizes among all jobs.\\
	$S:=\{S_1,\ldots, S_{\theta}\}$ : set of distinct sizes of jobs.\\
	$NT_{\ell,t}:=\displaystyle\sum_{\substack{ j \in J: s_j=S_\ell, \\ p_j=P_t}} 1$ : number of jobs of size $S_\ell$ and processing time $P_t$, for    $\ell=1,\ldots,\theta$, $t\in T$.\\
	$NT^+_{\ell,t}:= \displaystyle\sum_{\substack{ j \in J: s_j=S_\ell,\\ p_j\leq P_t}} 1$ : number of jobs of size $S_\ell$ and processing time $\leq P_t$,  for   $\ell=1,\ldots,\theta$, $t\in T$.\\
	$NJ_{t}:= \displaystyle\sum_{ \substack{j \in J:\\ p_j = P_t}} 1$ : number of jobs with processing time $P_t$,  for   $t\in T$.
%
%

\vspace{0.1in} 
\noindent
The set of nodes $V$ and the sets of arcs $A^J$ and $A^L$ are defined respectively in \eqref{defV}, \eqref{arc_J}, and \eqref{arc_L}.
	
\vspace{0.1in}

	\noindent
	\textbf{Decision variables}
	
	  	$f_{i,j,t}$ : flow on \emph{job arc} $(i,j) \in A^J$ in arc-flow structure $t$. The variable indicates the quantity of batches created with position $i$ occupied by jobs with size $j-i$. 
	
	$y_{i,j,t}$ : flow on the \emph{loss arc} $(i,B) \in A^L$ in arc-flow structure $t$.
	
    $v_{t}$ :  flow  on the \emph{feedback arc}  in arc-flow structure  $t$. The variable indicates the number of batches required with processing time $P_t$.
	
    $z_{c,t}$ : number of jobs with size $c$, not allocated in the batches with processing time smaller than or equal to $P_t$. Theses jobs are allowed to be allocated in the batches with processing time $P_{t+1}$.

\begin{fleqn}
\begin{align}
  (\mbox{FLOW})\;\;&\min \; \sum_{\forall t \in T} P_t.v_t & \label{ModelFlow_Cmax}\\
  \nonumber
    &\left(\sum_{\substack{(i,j) \in A^J}} f_{i,j,t} + \sum_{\substack{(i,j) \in A^L}} y_{i,j,t} \right) - \\
    &\left(\sum_{\substack{(j,i) \in A^J}} f_{j,i,t} + \sum_{\substack{(j,i) \in A^L}} y_{j,i,t} \right) = 
     \left\{ \begin{array}{ll}
     -v_t & \mbox{if } j=0;\\
      v_t & \mbox{if } j=B;\\
      0 & \mbox{if } 0<j<B.
     \end{array} \right.
     &   t \in T
     \label{ModelFlow_Flow} \\
    & NT_{c,t} - \sum_{\substack{(i,j) \in A^J:\\j-i=c}} f_{i,j,t} =  
    \left\{ \begin{array}{ll}
     z_{c,t} & \mbox{if $t=1$};\\
     - z_{c,t-1} & \mbox{if $t=\delta$};\\
     z_{c,t}-z_{c,t-1} & \mbox{if $1<t<\delta$}.
     \end{array} \right.  &   
      c \in \{1..B\} 
     \label{ModelFlow_Allocated}
\end{align}
\vspace{-7mm}
\begin{align}
    \;\;\;\;\;\;\;\;\;\;\;\;\;\;\;
    &f_{i,j,t} \leq min(NJ_{t}, NT^+_{j-i,t}), \, f_{i,j,t} \in \mathbb{Z}  &   t \in T, (i,j) \in A^J   
    \label{ModelFlow_Var_F}\\
&v_t \leq NJ_{t}, \, v_t \in \mathbb{Z} &    t \in T  \label{ModelFlow_Var_v}\\
&    y_{i,j,t} \leq NJ_{t}, \, y_{i,j,t} \in \mathbb{Z} &   t \in T, (i,j) \in A^L  
    \label{ModelFlow_Var_y}\\
    &z_{c,t}    \leq NT^+_{c,t}, \, z_{c,t}   \in \mathbb{Z} &  t \in T:t<\delta,   c \in \{1..B\} 
    \label{ModelFlow_Var_z}
\end{align}
\end{fleqn}
  
The objective function \eqref{ModelFlow_Cmax} minimizes the makespan.
Constraints \eqref{ModelFlow_Flow} ensure the flow conservation.
Constraints \eqref{ModelFlow_Allocated} ensure that all jobs are assigned and also control the number of jobs of each size that  are assigned to each arc-flow structure.
  Constraints \eqref{ModelFlow_Var_F}, \eqref{ModelFlow_Var_v}, \eqref{ModelFlow_Var_y} and \eqref{ModelFlow_Var_z} define the domains of the variables and their respective upper bounds, defined by Rules \ref{rulearc3}--\ref{rulearc4}.

\section{Computational results}
\label{ComputationalResults}

In this section, we discuss the numerical experiments where we compare the new arc-flow formulation (FLOW) for problem  $1|s_j,B|C_{\max}$, to the formulations from the literature,  (MILP$_1$)  and  (MILP$_1^+$). Three sets of test instances were considered in our experiments. 
The first and second sets were  respectively created by the authors of \cite{Chen2011} and \cite{Muter2020},  who kindly made them available.
The third set was created in this work, with instances randomly generated. The selection of the parameters for these instances aimed at  making them a bigger challenge for our formulation.
The experiments were executed on a computer with a 2.70GHz \texttt{Intel Quad-Core Xeon E5-2697 v2} processor and 64GB of RAM, using \texttt{CPLEX}, version 12.7.1.0. The computational time was limited to 1800 seconds for each instance.

\subsection{Instances from Chen et al. (2011) } 
  
 The first set of test instances for problem $1|s_j,B|C_{\max}$ is the same one considered in \cite{Chen2011}. 
To generate the instances used in our tests, six different numbers of jobs ($n_J$) were considered, as well as two intervals of integers  ($p_1$,$p_2$), from which  the processing times were randomly selected,  and three  intervals of integers ($s_1$,$s_2$,$s_3$), from which the sizes of the jobs were randomly selected. Considering the parameters shown in Table \ref{tabAllParameters}, instances of $36 (6 \times 2 \times 3)$ configurations were tested. For each configuration, 100 instances were generated. 
We solved the test instances with  \texttt{CPLEX}, configured to run in only one thread, not to benefit from the processor parallelism.

\begin{singlespace}
\begin{table}[!ht] \scalefont{1}
	\caption{Parameter setting for instances from \cite{Chen2011}.}
	\centering
	\begin{tabular}{ c c c c }
		\toprule
		Number of jobs   & Processing time & Job size     & Machine capacity \\ 
		($n_J$) & ($p_j$) & ($s_j$) & ($B$)\\ 
		\midrule
		10, 20, 50, & $p_1$: [1, 10]  &\multicolumn{1}{l}{$s_1$: [1, 10]} & $B=10$ \\
		100, 300, 500    & $p_2$: [1, 20]  & \multicolumn{1}{l}{$s_2$: [2, 4]}  &        \\
		                 &                 & \multicolumn{1}{l}{$s_3$: [4, 8]}  &        \\
		\bottomrule

	\end{tabular}
	\label{tabAllParameters}
\end{table}
\end{singlespace}

The following statistics were considered in our comparative  analysis: 
the computational time of \texttt{CPLEX} in seconds ($T(s)$) (we represent the time by the symbol ``-'' in our tables when \texttt{CPLEX} reaches the time limit of 1800 seconds on all instances of a given configuration), 
the  makespan corresponding to the best solution obtained by \texttt{CPLEX}  ($C_{\max}$),
the  duality gap of \texttt{CPLEX} at the end of its execution (Gap),
 and the  number of instances for which  \texttt{CPLEX} finds the optimal solution (\#O).

\begin{table}[!pt] 
\centering
\scalefont{1}
\setlength{\tabcolsep}{2.3pt}
\caption{Computational results for the instances available by \cite{Chen2011} for the $1|s_j.B|C_{\max}$ problem.}{
\begin{tabular}{ c c | r r r r | r r r r | r r r r }
\toprule
\multicolumn{2}{c}{Instance} & \multicolumn{4}{c}{(MILP$_1$)} & \multicolumn{4}{c}{(MILP$_1^+$)} & \multicolumn{4}{c}{(FLOW)} \\
\cmidrule(l){3-6} \cmidrule(l){7-10} \cmidrule(l){11-14}
\multicolumn{1}{c}{Jobs} & 
\multicolumn{1}{c}{Type} & 
\multicolumn{1}{c}{$T(s)$} & 
\multicolumn{1}{c}{$C_{\max}$} & 
\multicolumn{1}{c}{Gap} & 
\multicolumn{1}{c}{\#O} & 
\multicolumn{1}{c}{$T(s)$} & 
\multicolumn{1}{c}{$C_{\max}$} & 
\multicolumn{1}{c}{Gap} & 
\multicolumn{1}{c}{\#O} & 
\multicolumn{1}{c}{$T(s)$} & 
\multicolumn{1}{c}{$C_{\max}$} & 
\multicolumn{1}{c}{Gap} & 
\multicolumn{1}{c}{\#O} \\ 
\midrule
\multicolumn{14}{c}{Instances with $p_1 = [1.10]$}\\
\midrule
10 & $p_1s_1$  & 0.08 & \textbf{36.86}  & \textbf{0.00}  & \textbf{100}  & \textbf{0.00}  & \textbf{36.86}  & \textbf{0.00}  & \textbf{100}  & 0.01 & \textbf{36.86}  & \textbf{0.00}  & \textbf{100} \\
10 & $p_1s_2$  & 0.04 & \textbf{20.38}  & \textbf{0.00}  & \textbf{100}  & \textbf{0.01}  & \textbf{20.38}  & \textbf{0.00}  & \textbf{100}  & 0.02 & \textbf{20.38}  & \textbf{0.00}  & \textbf{100} \\
10 & $p_1s_3$  & 0.1 & \textbf{43.79}  & \textbf{0.00}  & \textbf{100}  & \textbf{0.00}  & \textbf{43.79}  & \textbf{0.00}  & \textbf{100}  & \textbf{0.00}  & \textbf{43.79}  & \textbf{0.00}  & \textbf{100} \\
20 & $p_1s_1$  & 42.62 & \textbf{68.07}  & 0.06 & 99 & \textbf{0.01}  & \textbf{68.07}  & \textbf{0.00}  & \textbf{100}  & 0.03 & \textbf{68.07}  & \textbf{0.00}  & \textbf{100} \\
20 & $p_1s_2$  & 300.49 & \textbf{37.13}  & 1.17 & 92 & 0.09 & \textbf{37.13}  & \textbf{0.00}  & \textbf{100}  & \textbf{0.06}  & \textbf{37.13}  & \textbf{0.00}  & \textbf{100} \\
20 & $p_1s_3$  & 24.37 & \textbf{83.69}  & 0.03 & 99 & \textbf{0.00}  & \textbf{83.69}  & \textbf{0.00}  & \textbf{100}  & 0.01 & \textbf{83.69}  & \textbf{0.00}  & \textbf{100} \\
50 & $p_1s_1$  & 1618.35 & 164.17 & 7.83 & 15 & 0.84 & \textbf{164.08}  & \textbf{0.00}  & \textbf{100}  & \textbf{0.18}  & \textbf{164.08}  & \textbf{0.00} & \textbf{100} \\
50 & $p_1s_2$  &       -  & 87.94 & 66.65 & 0 & 317.12 & \textbf{87.39}  & 0.2 & 88 & \textbf{0.36}  & \textbf{87.39}  & \textbf{0.00}  & \textbf{100} \\
50 & $p_1s_3$  & 1649.45 & \textbf{202.03}  & 14.17 & 13 & 0.02 & \textbf{202.03}  & \textbf{0.00}  & \textbf{100}  & \textbf{0.01}  & \textbf{202.03}  & \textbf{0.00}  & \textbf{100} \\
100 & $p_1s_1$  & -  & 325.38 & 73.99 & 0 & 70.24 & \textbf{318.99}  & 0.02 & 97 & \textbf{0.17}  & \textbf{318.99}  & \textbf{0.00}  & \textbf{100} \\
100 & $p_1s_2$  & -  & 183.78 & 88.82 & 0 & 1689.67 & \textbf{170.58}  & 1.41 & 7 & \textbf{0.61}  & \textbf{170.58}  & \textbf{0.00}  & \textbf{100} \\
100 & $p_1s_3$  & -  & 401.7 & 83.99 & 0 & 0.11 & \textbf{396.96}  & \textbf{0.00}  & \textbf{100}  & \textbf{0.01}  & \textbf{396.96}  & \textbf{0.00}  & \textbf{100} \\
300 & $p_1s_1$  & -  & 1832.31 & 99.38 & 0 & 467.18 & 928.64 & 0.05 & 76 & \textbf{0.23}  & \textbf{928.63}  & \textbf{0.00}  & \textbf{100} \\
300 & $p_1s_2$  & -  & 2508.79 & 99.53 & 0 &      -  & 496.07 & 0.86 & 0 & \textbf{1.04}  & \textbf{495.66}  & \textbf{0.00}  & \textbf{100} \\
300 & $p_1s_3$  & -  & 1999.85 & 99.44 & 0 & 20.4 & \textbf{1174.46}  & 0.0009 & 99 & \textbf{0.09}  & \textbf{1174.46}  & \textbf{0.00}  & \textbf{100} \\
500 & $p_1s_1$  & -  & 3112.85 & 99.77 & 0 & 825.1 & 1544.31 & 0.05 & 59 & \textbf{0.17}  & \textbf{1544.30}  & \textbf{0.00}  & \textbf{100} \\
500 & $p_1s_2$  & -  & 4970.07 & 99.8 & 0 &      -  & 832.43 & 0.72 & 0 & \textbf{0.97}  & \textbf{831.04}  & \textbf{0.00}  & \textbf{100} \\
500 & $p_1s_3$  & -  & 3347.23 & 99.88 & 0 & 34.39 & \textbf{1949.76}  & 0.0005 & 99 & \textbf{0.02}  & \textbf{1949.76}  & \textbf{0.00}  & \textbf{100} \\
\midrule
\multicolumn{14}{c}{Instances with $p_2 = [1.20]$}\\
\midrule
10 & $p_2s_1$  & 0.07 & \textbf{67.62}  & \textbf{0.00}  & \textbf{100}  & \textbf{0.00}  & \textbf{67.62}  & \textbf{0.00}  & \textbf{100}  & 0.02 & \textbf{67.62}  & \textbf{0.00}  & \textbf{100} \\
10 & $p_2s_2$  & 0.04 & \textbf{40.22}  & \textbf{0.00}  & \textbf{100}  & \textbf{0.02}  & \textbf{40.22}  & \textbf{0.00}  & \textbf{100}  & 0.03 & \textbf{40.22}  & \textbf{0.00}  & \textbf{100} \\
10 & $p_2s_3$  & 0.09 & \textbf{81.05}  & \textbf{0.00}  & \textbf{100}  & \textbf{0.00}  & \textbf{81.05}  & \textbf{0.00}  & \textbf{100}  & 0.01 & \textbf{81.05}  & \textbf{0.00}  & \textbf{100} \\
20 & $p_2s_1$  & 22.35 & \textbf{133.09}  & 0.01 & 99 & \textbf{0.01}  & \textbf{133.09}  & \textbf{0.00}  & \textbf{100}  & 0.08 & \textbf{133.09}  & \textbf{0.00}  & \textbf{100} \\
20 & $p_2s_2$  & 215.7 & \textbf{72.88}  & 0.96 & 94 & \textbf{0.09}  & \textbf{72.88}  & \textbf{0.00}  & \textbf{100}  & 0.13 & \textbf{72.88}  & \textbf{0.00}  & \textbf{100} \\
20 & $p_2s_3$  & 16.49 & \textbf{159.11}  & \textbf{0.00}  & \textbf{100}  & \textbf{0.00}  & \textbf{159.11}  & \textbf{0.00}  & \textbf{100}  & 0.01 & \textbf{159.11}  & \textbf{0.00}  & \textbf{100} \\
50 & $p_2s_1$  & 1596.95 & 314.76 & 9.4 & 15 & \textbf{0.38}  & \textbf{314.57}  & \textbf{0.00}  & \textbf{100}  & 0.62 & \textbf{314.57}  & \textbf{0.00}  & \textbf{100} \\
50 & $p_2s_2$  &       -  & 169.93 & 65.57 & 0 & 178.7 & \textbf{168.11}  & 0.06 & 94 & \textbf{1.81}  & \textbf{168.11}  & \textbf{0.00}  & \textbf{100} \\
50 & $p_2s_3$  & 1701.53 & 384.15 & 15.56 & 8 & \textbf{0.02}  & \textbf{384.13 } & \textbf{0.00}  & \textbf{100}  & \textbf{0.02}  & \textbf{384.13}  & \textbf{0.00}  & \textbf{100} \\
100 & $p_2s_1$  & -  & 622.84 & 76.62 & 0 & 33.35 & \textbf{610.64}  & 0.01 & 99 & \textbf{0.86}  & \textbf{610.64}  & \textbf{0.00}  & \textbf{100} \\
100 & $p_2s_2$  & -  & 357.69 & 89.35 & 0 & 1717.84 & 326.14 & 0.99 & 6 & \textbf{4.71}  & \textbf{326.11}  & \textbf{0.00}  & \textbf{100} \\
100 & $p_2s_3$  & -  & 775.5 & 83.96 & 0 & 0.1 & \textbf{766.91}  & \textbf{0.00}  & \textbf{100}  & \textbf{0.04}  & \textbf{766.91}  & \textbf{0.00}  & \textbf{100} \\
300 & $p_2s_1$  & -  & 3979.09 & 99.42 & 0 & 329.03 & 1793.54 & 0.02 & 83 & \textbf{0.96}  & \textbf{1793.52}  & \textbf{0.00}  & \textbf{100} \\
300 & $p_2s_2$  & -  & 5925.99 & 99.66 & 0 &      -  & 965.42 & 1.08 & 0 & \textbf{6.55}  & \textbf{962.77}  & \textbf{0.00}  & \textbf{100} \\
300 & $p_2s_3$  & -  & 3983.88 & 99.43 & 0 & 21.78 & \textbf{2247.39}  & 0.0005 & 99 & \textbf{0.17}  & \textbf{2247.39}  & \textbf{0.00}  & \textbf{100} \\
500 & $p_2s_1$  & -  & 6264.18 & 99.92 & 0 & 780.62 & 2964.62 & 0.03 & 63 & \textbf{0.89}  & \textbf{2964.57}  & \textbf{0.00}  & \textbf{100} \\
500 & $p_2s_2$  & -  & 9937.42 & 99.94 & 0 &      -  & 1592.72 & 0.86 & 0 & \textbf{5.66}  & \textbf{1587.50}  & \textbf{0.00}  & \textbf{100} \\
500 & $p_2s_3$  & -  & 6822.02 & 99.92 & 0 & 97.16 & \textbf{3701.79}  & 0.0013 & 96 & \textbf{0.13}  & \textbf{3701.79}  & \textbf{0.00}  & \textbf{100} \\
\toprule
\end{tabular}
}
\label{Cap3tab10jobs}
\end{table}
  
We present in Table \ref{Cap3tab10jobs} comparison results among the three formulations for problem  $1|s_j,B|C_{\max}$. 
All values presented in  Table \ref{Cap3tab10jobs} are average results computed over the instances of the same configuration, as described in Table \ref{tabAllParameters}. 
We note that it is possible to have the computational times of \texttt{CPLEX}, for solving problems of a given configuration, less than 1800 seconds while the gaps are non-zero. This happens when some of the instances in the group could be solved to optimality in the time limit, and others could not. In this case, the average considers all the instances in the group, i.e., we sum up, the time of 1800 seconds for the instances that could not be solved to optimality and the smaller times for instances that could, and divide the sum by the total number of instances in the group.

The comparative tests clearly show that formulation (FLOW) is far superior to (MILP$_1$) and (MILP$_1^+$) for these instances,  especially when the number of jobs increases.
This is the first time optimal solutions are shown for all these instances, with up to 500 jobs.
Considering the small instances, with 50 jobs or less, we see that  (MILP$_1^+$) can solve some of the instances in less computational time than (FLOW), but the difference between times is always a fraction of a second.
Additionally, the duality gaps shown for (MILP$_1$) reveal the difficulty in obtaining good lower bounds with this formulation. 
This difficulty is reduced with the use of (MILP$_1^+$), but only (FLOW) can show results with no gap for any instances with 50 jobs or more.
With formulation (FLOW), we were able to prove the optimality for 86\% of the instances in less than 1 second, while with models (MILP$_1^+$) and (MILP$_1$), we  prove optimality for 56.75\% and 17.28\% of the instances, respectively.

The total computational time spent on our 3600 test instances, when using model (FLOW), was 44 minutes and 29 seconds,  while it was 15 days, 22 hours, 53 minutes and 56 seconds when using (MILP$_1^+$), and 49 days, 23 hours, 41 minutes and 13 seconds, when using (MILP$_1$).

Unlike what occurs with models (MILP$_1$) and (MILP$_1^+$), the number of variables in (FLOW) does not increase when the number of jobs increases. Moreover, the flow graph does not change in this case. Only the bounds on the variables change. The flow graphs of two distinct instances will be the same if the settings in the parameters  Processing Time, Job Size, and Machine Capacity are the same. In fact, this is a very important characteristic of the flow approach. We finally note that the computational time to construct the graphs for the flow formulation was not considered in these times. However, the maximum time to construct a graph for any instance in these experiments was 0.008 second.

Table \ref{Table_Results1_Relax} shows the solutions of the linear relaxations of models (MILP$_1^+$) and (FLOW).
The following statistics were considered: 
the makespan corresponding to the optimal solution of the instance previous calculated ($C_{\max}*$),
the makespan corresponding to the solution of the relaxations ($C_{\max}$),
the gap given by: gap = ($C_{\max}*$ - $C_{\max}$)/$C_{\max}$,
the number of simplex iterations (Iter.), and
the computational time in seconds ($T(s)$).
\begin{table}[!p] 
\centering
\scalefont{1}
\setlength{\tabcolsep}{5.5pt}
\caption{Comparison between the linear relaxations for (MILP$_1^+$) and (FLOW) for the $1|s_j,B|C_{\max}$ problem.}{
\begin{tabular}{ c c r | r r r r | r r r r }
\toprule
\multicolumn{3}{c}{Instance} & 
\multicolumn{4}{c}{(MILP$_1^+$) - Relax} & 
\multicolumn{4}{c}{(FLOW) - Relax} \\
\cmidrule(l){4-7} \cmidrule(l){8-11}
\multicolumn{1}{c}{Jobs} & 
\multicolumn{1}{c}{Type} & 
\multicolumn{1}{c}{$C_{max}*$} & 
\multicolumn{1}{c}{$C_{\max}$} & 
\multicolumn{1}{c}{Gap} & 
\multicolumn{1}{c}{Iter.} & 
\multicolumn{1}{c}{$T(s)$} & 
\multicolumn{1}{c}{$C_{\max}$} & 
\multicolumn{1}{c}{Gap} & 
\multicolumn{1}{c}{Iter.} & 
\multicolumn{1}{c}{$T(s)$} \\
\midrule
\multicolumn{11}{c}{Instances with $p_1 = [1,10]$}\\
\midrule
10 & $p_1s_1$  & 36.86 & 32.34 & 12.26 & \textbf{29.41} & \textbf{0.00} & \textbf{34.41} & \textbf{6.66} & 79.52 & \textbf{0.00}\\
10 & $p_1s_2$  & 20.38 & \textbf{18.61} & \textbf{8.70} & \textbf{36.66} & \textbf{0.00} & 16.79 & 17.61 & 82.88 & \textbf{0.00}\\
10 & $p_1s_3$  & 43.79 & 34.82 & 20.49 & 34.09 & \textbf{0.00} & \textbf{41.72} & \textbf{4.74} & \textbf{14.77} & \textbf{0.00}\\
20 & $p_1s_1$  & 68.07 & 61.42 & 9.77 & \textbf{80.21} & \textbf{0.00} & \textbf{65.34} & \textbf{4.02} & 181.54 & \textbf{0.00}\\
20 & $p_1s_2$  & 37.13 & \textbf{34.54} & \textbf{6.99} & \textbf{103.60} & \textbf{0.00} & 33.38 & 10.11 & 133.90 & \textbf{0.00}\\
20 & $p_1s_3$  & 83.69 & 66.80 & 20.18 & 96.61 & \textbf{0.00} & \textbf{81.49} & \textbf{2.63} & \textbf{26.04} & \textbf{0.00}\\
50 & $p_1s_1$  & 164.08 & 151.94 & 7.40 & 330.65 & 0.01 & \textbf{161.03} & \textbf{1.86} & \textbf{318.89} & \textbf{0.01}\\
50 & $p_1s_2$  & 87.39 & 83.07 & 4.94 & 341.59 & 0.01 & \textbf{83.15} & \textbf{4.85} & \textbf{168.99} & \textbf{0.00}\\
50 & $p_1s_3$  & 202.03 & 165.65 & 18.01 & 385.43 & 0.01 & \textbf{199.71} & \textbf{1.15} & \textbf{42.95} & \textbf{0.00}\\
100 & $p_1s_1$  & 318.99 & 302.63 & 5.13 & 675.16 & 0.04 & \textbf{315.98} & \textbf{0.94} & \textbf{386.88} & \textbf{0.01}\\
100 & $p_1s_2$  & 170.58 & 165.68 & 2.87 & 720.13 & 0.05 & \textbf{165.99} & \textbf{2.69} & \textbf{178.51} & \textbf{0.00}\\
100 & $p_1s_3$  & 396.96 & 328.82 & 17.17 & 742.44 & 0.05 & \textbf{394.65} & \textbf{0.58} & \textbf{56.13} & \textbf{0.00}\\
300 & $p_1s_1$  & 928.63 & 902.86 & 2.77 & 2729.99 & 0.36 & \textbf{925.87} & \textbf{0.30} & \textbf{467.62} & \textbf{0.01}\\
300 & $p_1s_2$  & 495.66 & 490.81 & 0.98 & 3186.63 & 0.44 & \textbf{490.85} & \textbf{0.97} & \textbf{180.92} & \textbf{0.01}\\
300 & $p_1s_3$  & 1174.46 & 988.13 & 15.87 & 3043.30 & 0.40 & \textbf{1172.12} & \textbf{0.20} & \textbf{62.20} & \textbf{0.00}\\
500 & $p_1s_1$  & 1544.30 & 1513.43 & 2.00 & 4673.94 & 1.19 & \textbf{1541.61} & \textbf{0.17} & \textbf{480.35} & \textbf{0.01}\\
500 & $p_1s_2$  & 831.04 & 826.12 & 0.59 & 5387.19 & 1.64 & \textbf{826.12} & \textbf{0.59} & \textbf{180.45} & \textbf{0.01}\\
500 & $p_1s_3$  & 1949.76 & 1645.82 & 15.59 & 5354.30 & 1.47 & \textbf{1947.47} & \textbf{0.12} & \textbf{61.54} & \textbf{0.00}\\
\midrule
\multicolumn{11}{c}{Instances with $p_2 = [1,20]$}\\
\midrule
10 & $p_2s_1$  & 67.62 & 59.67 & 11.75 & \textbf{30.16} & \textbf{0.00} & \textbf{62.33} & \textbf{7.82} & 109.14 & \textbf{0.00}\\
10 & $p_2s_2$  & 40.22 & \textbf{36.55} & \textbf{9.14} & \textbf{37.14} & \textbf{0.00} & 32.69 & 18.72 & 101.74 & \textbf{0.00}\\
10 & $p_2s_3$  & 81.05 & 64.90 & 19.93 & 33.71 & \textbf{0.00} & \textbf{77.00} & \textbf{5.00} & \textbf{18.71} & \textbf{0.00}\\
20 & $p_2s_1$  & 133.09 & 119.97 & 9.86 & \textbf{84.59} & \textbf{0.00} & \textbf{128.09} & \textbf{3.75} & 253.47 & 0.01\\
20 & $p_2s_2$  & 72.88 & \textbf{68.24} & \textbf{6.37} & \textbf{105.44} & \textbf{0.00} & 65.23 & 10.50 & 181.63 & \textbf{0.00}\\
20 & $p_2s_3$  & 159.11 & 129.26 & 18.76 & 97.66 & \textbf{0.00} & \textbf{154.68} & \textbf{2.78} & \textbf{38.11} & \textbf{0.00}\\
50 & $p_2s_1$  & 314.57 & 291.85 & 7.22 & \textbf{331.51} & \textbf{0.01} & \textbf{309.22} & \textbf{1.70} & 566.30 & 0.02\\
50 & $p_2s_2$  & 168.11 & \textbf{161.19} & \textbf{4.12} & 352.33 & 0.01 & 159.89 & 4.89 & \textbf{295.78} & \textbf{0.00}\\
50 & $p_2s_3$  & 384.13 & 313.80 & 18.31 & 374.71 & 0.01 & \textbf{379.53} & \textbf{1.20} & \textbf{73.38} & \textbf{0.00}\\
100 & $p_2s_1$  & 610.64 & 579.74 & 5.06 & \textbf{693.34} & 0.05 & \textbf{605.19} & \textbf{0.89} & 782.51 & \textbf{0.02}\\
100 & $p_2s_2$  & 326.11 & 316.77 & 2.86 & 742.66 & 0.05 & \textbf{317.20} & \textbf{2.73} & \textbf{335.60} & \textbf{0.01}\\
100 & $p_2s_3$  & 766.91 & 633.63 & 17.38 & 783.24 & 0.05 & \textbf{762.26} & \textbf{0.61} & \textbf{102.42} & \textbf{0.00}\\
300 & $p_2s_1$  & 1793.52 & 1737.65 & 3.11 & 2786.36 & 0.36 & \textbf{1788.23} & \textbf{0.29} & \textbf{963.18} & \textbf{0.03}\\
300 & $p_2s_2$  & 962.77 & 952.38 & 1.08 & 3181.75 & 0.41 & \textbf{952.79} & \textbf{1.04} & \textbf{357.50} & \textbf{0.01}\\
300 & $p_2s_3$  & 2247.39 & 1892.29 & 15.80 & 3161.54 & 0.41 & \textbf{2242.72} & \textbf{0.21} & \textbf{135.67} & \textbf{0.00}\\
500 & $p_2s_1$  & 2964.57 & 2895.60 & 2.33 & 4966.34 & 1.25 & \textbf{2959.25} & \textbf{0.18} & \textbf{1017.03} & \textbf{0.04}\\
500 & $p_2s_2$  & 1587.50 & 1577.43 & 0.63 & 5496.04 & 1.51 & \textbf{1577.56} & \textbf{0.63} & \textbf{360.63} & \textbf{0.01}\\
500 & $p_2s_3$  & 3701.79 & 3138.60 & 15.21 & 5664.89 & 1.52 & \textbf{3697.31} & \textbf{0.12} & \textbf{136.55} & \textbf{0.00}\\
\toprule
\end{tabular}
}
\label{Table_Results1_Relax}
\end{table}

We notice that the linear relaxation of (FLOW) is much better than the linear relaxation of (MILP$_1^+$) for most instances, especially when the number of jobs increases.
For instances of type $s_2$, the solutions of the linear relaxations are very close to each other, but (FLOW) is slightly better in this case.
Another critical point is that (FLOW) presents lower computational time for solving the linear relaxations. 
For instances with 500 jobs and  type $p_2s_2$, the time for linear relaxation of (FLOW) is 0.66\% of the time for the linear relaxation of (MILP$_1^+$). 
This time difference tends to increase when the number of jobs increases.
These results give us some insight into why (FLOW) overcomes (MILP$_1^+$).

The results in Tables  \ref{Cap3tab10jobs}  and  \ref{Table_Results1_Relax} show that instances of configuration $s_2$ require more computational time and are more difficult when compared to the other instances for all formulations.
The reason for this is the small sizes of the jobs when compared to the machine capacity, which allows more combinations of assignment to a batch. In this way, the number of feasible solutions for instances of type $s_2$ is greater than for the other instances.

With the new formulation  (FLOW), it is possible to find the optimal solution of all the instances proposed by \cite{Chen2011} with good computational times. 

To illustrate the application of our approach to instances with a large number of jobs, where multiple jobs have equal size and processing time,  
we created new instances where the number of jobs goes up to 100 million, but the the intervals from which the processing times and sizes are selected, were maintained respectively, as $p_2$ and $s_2$. We aim at emphasizing here how the formulation is particularly suitable for this type of instance.
Table \ref{Table_Results1_Lit} shows the computational results obtained. 
With (FLOW), we could find the optimum solution of all the  instances with up to 100 million jobs in no more than 6.05 seconds.  
We  note that column ``Construction Time'' in Table \ref{Table_Results1_Lit} considers the times to  read the data of the instance,  sort the jobs and  create the parameters to the model, i.e., to construct the graph $G$. 
The algorithm to construct $G$ is  polynomial in the number of jobs, but the construction phase can take longer than the resolution of the problem, when the number of jobs increases.
Although the construction time increases with the number of jobs, even when added to the solution time of  model (FLOW), it leads to very small times to solve instances of this particular type, when comparing to other formulations from the literature.
The column ``Nodes'' show the number of nodes required by \texttt{CPLEX} to solve the problem.

\begin{table}[ht!] 
\centering
\scalefont{1}
\caption{Computational results for large random instances for the $1|s_j,B|C_{\max}$ problem.}{
\begin{tabular}{ r c | r r r r r }
\toprule
\multicolumn{2}{c}{Instance} & \multicolumn{5}{c}{(FLOW)} \\
\cmidrule(l){3-7} 
\multicolumn{1}{c}{Jobs} & 
\multicolumn{1}{c}{Type} & 
\multicolumn{1}{c}{$T(s)$} & 
\multicolumn{1}{c}{$C_{\max}$} & 
\multicolumn{1}{c}{Gap} & 
\multicolumn{1}{c}{Nodes} & 
\multicolumn{1}{c}{Construction Time}\\
\midrule
1000 & $p_2s_2$ & 4.90 & 3241 & 0.00 & 8043 & 0.03 \\
10000 & $p_2s_2$ & 3.66 & 31719 & 0.00 & 3625 & 0.16 \\
100000 & $p_2s_2$ & 3.22 & 314945 & 0.00 & 4668 & 1.44 \\
1000000 & $p_2s_2$ & 2.80 & 3152697 & 0.00 & 3011 & 14.24 \\
10000000 & $p_2s_2$ & 6.05 & 31495193 & 0.00 & 17111 & 145.07 \\
100000000 & $p_2s_2$ & 3.25 & 314996812 & 0.00 & 4325 & 886.14 \\
\toprule
\end{tabular}
}
\label{Table_Results1_Lit}
\end{table}

\subsection{Instances from Muter (2020) } 
\label{sec:InstancesSingle2}

The second set of test instances for problem $1|s_j,B|C_{\max}$ is the same one considered in the recent work \cite{Muter2020}. 
To generate these instances, four different numbers of jobs ($n_J$) were considered, as well as three intervals of integers  ($p_1$,$p_2$,$p_3$), from which  the processing times were randomly selected,  and three  intervals of integers ($s_1$,$s_2$,$s_3$), from which the sizes of the jobs were randomly selected. Considering the parameters shown in Table \ref{tabAllParametersM}, instances of $36(4 \times 3 \times 3)$ configurations were tested. For each configuration, 20 instances were generated. 
We solved the  instances with  \texttt{CPLEX} configured in the default settings. We emphasize that we use in these tests similar configuration to the one  used in \cite{Muter2020}.

\begin{singlespace}
\begin{table}[!ht] \scalefont{1}
	\caption{Parameter setting for instances from \cite{Muter2020}.}
	\centering
	\begin{tabular}{ c c c c }
		\toprule
		Number of jobs   & Processing time & Job size     & Machine capacity \\ 
		($n_J$) & ($p_j$) & ($s_j$) & ($B$)\\ 
		\midrule
		10, 20, 50, & $p_1$: [1, 10]  &\multicolumn{1}{l}{$s_1$: [1, 10]} & $B=10$ \\
		100         & $p_2$: [1, 20]  & \multicolumn{1}{l}{$s_2$: [2, 4]}  &        \\
		            & $p_3$: [1, 100] & \multicolumn{1}{l}{$s_3$: [4, 8]}  &        \\
		\bottomrule

	\end{tabular}
	\label{tabAllParametersM}
\end{table}
\end{singlespace}

We keep the notation used by the author to represent the instances in Table \ref{tableMuter}, where $n_j$ is the number of jobs, and $\bar{s}$ and $\bar{p}$ are the maximum sizes and processing times of the jobs in the interval.
The same statistics presented in Table \ref{Cap3tab10jobs} are presented  in Table \ref{tableMuter} comparing (MILP$_1^+$) and (FLOW) on this other set of instances. 

\begin{table}[!pt] 
\centering
\scalefont{1}
\setlength{\tabcolsep}{2.3pt}
\caption{Computational results for the instances available by \cite{Muter2020} for the $1|s_j.B|C_{\max}$ problem.}{
\begin{tabular}{ c c c | r r r r | r r r r }
\toprule
\multicolumn{3}{c}{Instance} & \multicolumn{4}{c}{(MILP$_1^+$)} & \multicolumn{4}{c}{(FLOW)} \\
\cmidrule(l){4-7} \cmidrule(l){8-11}
\multicolumn{1}{c}{$n_j$} & 
\multicolumn{1}{c}{$\bar{s}$} & 
\multicolumn{1}{c}{$\bar{p}$} & 
\multicolumn{1}{c}{$T(s)$} & 
\multicolumn{1}{c}{$C_{\max}$} & 
\multicolumn{1}{c}{Gap} & 
\multicolumn{1}{c}{\#O} & 
\multicolumn{1}{c}{$T(s)$} & 
\multicolumn{1}{c}{$C_{\max}$} & 
\multicolumn{1}{c}{Gap} & 
\multicolumn{1}{c}{\#O}  \\ 
\midrule
    10    & 4     & 10    & \textbf{0,05} & 21,45 & 0,00  & 20    & 0,08  & 21,45 & 0,00  & 20 \\
    10    & 4     & 20    & \textbf{0,05} & 39,15 & 0,00  & 20    & 0,10  & 39,15 & 0,00  & 20 \\
    10    & 4     & 100   & \textbf{0,06} & 202,35 & 0,00  & 20    & 0,11  & 202,35 & 0,00  & 20 \\
    10    & 8     & 10    & \textbf{0,01} & 44,85 & 0,00  & 20    & 0,02  & 44,85 & 0,00  & 20 \\
    10    & 8     & 20    & \textbf{0,02} & 78,25 & 0,00  & 20    & 0,03  & 78,25 & 0,00  & 20 \\
    10    & 8     & 100   & \textbf{0,01} & 414,10 & 0,00  & 20    & 0,03  & 414,10 & 0,00  & 20 \\
    10    & 10    & 10    & \textbf{0,03} & 34,25 & 0,00  & 20    & 0,05  & 34,25 & 0,00  & 20 \\
    10    & 10    & 20    & \textbf{0,02} & 71,35 & 0,00  & 20    & 0,05  & 71,35 & 0,00  & 20 \\
    10    & 10    & 100   & \textbf{0,02} & 333,85 & 0,00  & 20    & 0,06  & 333,85 & 0,00  & 20 \\
    20    & 4     & 10    & 0,10  & 37,35 & 0,00  & 20    & \textbf{0,09} & 37,35 & 0,00  & 20 \\
    20    & 4     & 20    & \textbf{0,12} & 72,90 & 0,00  & 20    & 0,16  & 72,90 & 0,00  & 20 \\
    20    & 4     & 100   & \textbf{0,13} & 357,95 & 0,00  & 20    & 0,34  & 357,95 & 0,00  & 20 \\
    20    & 8     & 10    & \textbf{0,02} & 82,65 & 0,00  & 20    & 0,04  & 82,65 & 0,00  & 20 \\
    20    & 8     & 20    & \textbf{0,02} & 159,05 & 0,00  & 20    & 0,04  & 159,05 & 0,00  & 20 \\
    20    & 8     & 100   & \textbf{0,04} & 811,65 & 0,00  & 20    & 0,06  & 811,65 & 0,00  & 20 \\
    20    & 10    & 10    & \textbf{0,04} & 66,85 & 0,00  & 20    & 0,07  & 66,85 & 0,00  & 20 \\
    20    & 10    & 20    & \textbf{0,04} & 123,45 & 0,00  & 20    & 0,10  & 123,45 & 0,00  & 20 \\
    20    & 10    & 100   & \textbf{0,06} & 639,95 & 0,00  & 20    & 0,21  & 639,95 & 0,00  & 20 \\
    50    & 4     & 10    & 55,52 & 88,30 & 0,00  & 20    & \textbf{0,31} & 88,30 & 0,00  & 20 \\
    50    & 4     & 20    & 13,30 & 171,95 & 0,00  & 20    & \textbf{0,96} & 171,95 & 0,00  & 20 \\
    50    & 4     & 100   & 20,39 & 819,15 & 0,00  & 20    & \textbf{6,33} & 819,15 & 0,00  & 20 \\
    50    & 8     & 10    & 0,07  & 205,05 & 0,00  & 20    & \textbf{0,05} & 205,05 & 0,00  & 20 \\
    50    & 8     & 20    & \textbf{0,04} & 397,80 & 0,00  & 20    & 0,06  & 397,80 & 0,00  & 20 \\
    50    & 8     & 100   & \textbf{0,06} & 1839,65 & 0,00  & 20    & 0,08  & 1839,65 & 0,00  & 20 \\
    50    & 10    & 10    & 0,15  & 163,70 & 0,00  & 20    & \textbf{0,11} & 163,70 & 0,00  & 20 \\
    50    & 10    & 20    & \textbf{0,17} & 316,50 & 0,00  & 20    & 0,23  & 316,50 & 0,00  & 20 \\
    50    & 10    & 100   & \textbf{0,19} & 1580,25 & 0,00  & 20    & 1,94  & 1580,25 & 0,00  & 20 \\
    100   & 4     & 10    & 1688,21 & 172,50 & 1,04  & 3     & \textbf{0,33} & 172,50 & \textbf{0,00} & \textbf{20} \\
    100   & 4     & 20    & 1670,83 & 319,50 & 0,90  & 2     & \textbf{2,49} & 319,50 & \textbf{0,00} & \textbf{20} \\
    100   & 4     & 100   & 928,37 & 1555,40 & 0,29  & 11    & \textbf{499,95} & \textbf{1555,35} & \textbf{0,07} & \textbf{17} \\
    100   & 8     & 10    & 0,15  & 399,35 & 0,00  & 20    & \textbf{0,06} & 399,35 & 0,00  & 20 \\
    100   & 8     & 20    & 0,17  & 766,75 & 0,00  & 20    & \textbf{0,07} & 766,75 & 0,00  & 20 \\
    100   & 8     & 100   & 0,16  & 3808,15 & 0,00  & 20    & \textbf{0,15} & 3808,15 & 0,00  & 20 \\
    100   & 10    & 10    & 1,71  & 312,50 & 0,00  & 20    & \textbf{0,18} & 312,50 & 0,00  & 20 \\
    100   & 10    & 20    & 0,97  & 597,85 & 0,00  & 20    & \textbf{0,29} & 597,85 & 0,00  & 20 \\
    100   & 10    & 100   & \textbf{0,86} & 2989,95 & 0,00  & 20    & 3,41  & 2989,95 & 0,00  & 20 \\
\toprule
\end{tabular}
}
\label{tableMuter}
\end{table}

For the second set of instances we have a similar analysis to the one obtained for the first one. 
The results  show that formulation (FLOW) is far superior to (MILP$_1^+$) on the most difficult instances, where the computation times for  (MILP$_1^+$) were more than  2 seconds.
Considering the small instances, with 20 jobs or less, we see that  (MILP$_1^+$) can solve almost all instances in less computational time than (FLOW), but the difference between times is a fraction of a second. Considering instances with 50 or 100 jobs, we see that the times used by (FLOW) is on average, smaller than the times used by (MILP$_1^+$). 
The  average time to solve all instances was 121.73 seconds for (MILP$_1^+$), compared to only 14.41 seconds for (FLOW).
However, we should note that the difficulty for (FLOW) increases when the processing time increases, which  is expected.  Concerning the duality gaps, the average values for (MILP$_1^+$) and (FLOW) are, respectively, 
0.06 and 0.002.
With formulation (FLOW), we were able to prove the optimality for 717 instances, i.e, for  99.6\% of the instances, while with (MILP$_1^+$) we  prove optimality for 676 instances,  93.9\% of them.

In \cite{Muter2020}, a column-and-cut generation CCG algorithm is applied to all the instances of the test set considered in Table \ref{tableMuter} that could not be solved to optimality in 300 seconds when using formulation (MILP$_1^+$). 
The results presented  in the paper show that the total number of instances solved increase from 675 to 694, still less than the 717 instances solved by (FLOW). 
The experiments reported in \cite{Muter2020} were conducted on the nodes of a computing cluster, each with two Intel E5-
2650 v2 CPUs (16 cores in total) running at 2.60 GHz. The CCG algorithm was also solved under a time limit of 1800 seconds.  

\subsection{New instances proposed } \label{sec:NewInstances}

From the computational tests with the instances  proposed by \cite{Chen2011}, we notice that these instances are not challenging enough for (FLOW). Also, from  the tests with the instances proposed by \cite{Muter2020},  we can see that high processing times increase the difficulty for the (FLOW) model.

Therefore, a third set of test instances for problem $1|s_j,B|C_{\max}$ is proposed in this work, with the parameter setting depicted in Table \ref{tabAllParameters1_2}.
The main idea is to generate more difficult instances, specifically for (FLOW), in order to better verify for which instances the model is more effective. For this, we explore the increase in the value of the parameter $B$, and in the size of the intervals from which $p_j$ and $s_j$ are selected. These increases  directly affect the number of variables and constraints  in formulation (FLOW).
Comparing with the set of instances from \cite{Chen2011} this new set of instances considers three larger values for the machine capacity $B$. The processing times are randomly selected from two intervals of integers again, but the second interval is much larger when we have a large  number of jobs. This will lead to more distinct values of processing times among the jobs, and therefore, to the necessity of more arc-flow structures to represent the solutions. Finally, job sizes  are randomly selected from three intervals of integers,  where the  ranges are proportional to the size of the machine, for example,   the interval $s_2 :=[ 0.2B ,  0.4B ]$. Note that when $B=100$, the interval is [$ 20 $, $ 40 $], so it keeps the  similar  characteristics of the instances from \cite{Chen2011}, in terms of scale.

In total, 540 instances  were generated, 5 for each of the 108 different combinations of number of jobs, range of processing times, range of job sizes, and machine capacity. 
We solved the test instances with  \texttt{CPLEX}, configured to run in only one thread, not to benefit from the processor parallelism.

\begin{table}[h] 
	\caption{Parameter setting for the new set of instances.}
	\centering
	\begin{tabular}{ c c c c }
		\toprule
		Number of jobs   & Processing time & Jobs size     & Machine capacity \\ 
		($n_J$) & ($p_j$) & ($s_j$) & ($B$)\\ 
		\midrule
		10, 50, 100, 500,& $p_1$: [1, 20]     & \multicolumn{1}{l}{$s_1$: [1, $B$]} & $20,50,100$ \\
		1000, 5000       & $p_2$: [1, $n_J$]  & \multicolumn{1}{l}{$s_2$: [$0.2B$, $0.4B$]}  &        \\
		                 &                    & \multicolumn{1}{l}{$s_3$: [$0.4B$, $0.8B$]}  &        \\
		\bottomrule

	\end{tabular}
	\label{tabAllParameters1_2}
\end{table}

We present in Tables \ref{Table_Results1_New_1}--\ref{Table_Results1_New_3} comparison results between formulations  (MILP$_1^+$) and (FLOW) for the new set of instances generated. When \texttt{CPLEX} cannot find an integer solution of a given configuration, we represent the $C_{\max}$ by ``No solution" and the gap by ``$\infty$" in our tables.

\begin{table}[!pt] 
\centering
\scalefont{1}
\setlength{\tabcolsep}{2.3pt}
\caption{Computational results for new instances proposed ($B=20$).}{
\begin{tabular}{ c c c | r r r r | r r r r }
\toprule
\multicolumn{3}{c}{Instance} & \multicolumn{4}{c}{(MILP$_1^+$)} & \multicolumn{4}{c}{(FLOW)} \\
\cmidrule(l){4-7} \cmidrule(l){8-11}
\multicolumn{1}{c}{Jobs} & 
\multicolumn{1}{c}{Capacity} & 
\multicolumn{1}{c}{Type} &
\multicolumn{1}{c}{$T(s)$} &
\multicolumn{1}{c}{$C_{\max}$} &
\multicolumn{1}{c}{Gap} & 
\multicolumn{1}{c}{\#O} &
\multicolumn{1}{c}{$T(s)$} & 
\multicolumn{1}{c}{$C_{\max}$} & 
\multicolumn{1}{c}{Gap} &
\multicolumn{1}{c}{\#O} \\ 
\midrule
\multicolumn{11}{c}{Instances with $p_1 = [1,20]$ and $B = 20$}\\
\midrule
10 & 20 & $p_1s_1$ & \textbf{0.01} & \textbf{62.20} & \textbf{0.00} & \textbf{5} & 0.03 & \textbf{62.20} & \textbf{0.00} & \textbf{5}\\
10 & 20 & $p_1s_2$ & \textbf{0.01} & \textbf{45.40} & \textbf{0.00} & \textbf{5} & 0.12 & \textbf{45.40} & \textbf{0.00} & \textbf{5}\\
10 & 20 & $p_1s_3$ & \textbf{0.00} & \textbf{71.80} & \textbf{0.00} & \textbf{5} & 0.01 & \textbf{71.80} & \textbf{0.00} & \textbf{5}\\
50 & 20 & $p_1s_1$ & \textbf{0.63} & \textbf{316.20} & \textbf{0.00} & \textbf{5} & 1.23 & \textbf{316.20} & \textbf{0.00} & \textbf{5}\\
50 & 20 & $p_1s_2$ & 276.05 & \textbf{181.00} & \textbf{0.00} & \textbf{5} & \textbf{13.67} & \textbf{181.00} & \textbf{0.00} & \textbf{5}\\
50 & 20 & $p_1s_3$ & \textbf{0.01} & \textbf{373.80} & \textbf{0.00} & \textbf{5} & 0.02 & \textbf{373.80} & \textbf{0.00} & \textbf{5}\\
100 & 20 & $p_1s_1$ & \textbf{0.87} & \textbf{629.60} & \textbf{0.00} & \textbf{5} & 1.20 & \textbf{629.60} & \textbf{0.00} & \textbf{5}\\
100 & 20 & $p_1s_2$ & 1468.32 & 326.80 & 0.90 & 2 & \textbf{20.91} & \textbf{326.40} & \textbf{0.00} & \textbf{5}\\
100 & 20 & $p_1s_3$ & 0.05 & \textbf{791.00} & \textbf{0.00} & \textbf{5} & \textbf{0.02} & \textbf{791.00} & \textbf{0.00} & \textbf{5}\\
500 & 20 & $p_1s_1$ & 1145.40 & 2805.40 & 0.03 & 2 & \textbf{2.52} & \textbf{2805.20} & \textbf{0.00} & \textbf{5}\\
500 & 20 & $p_1s_2$ & - & 1613.60 & 1.80 & 0 & \textbf{39.35} & \textbf{1595.20} & \textbf{0.00} & \textbf{5}\\
500 & 20 & $p_1s_3$ & 2.48 & \textbf{3869.80} & \textbf{0.00} & \textbf{5} & \textbf{0.07} & \textbf{3869.80} & \textbf{0.00} & \textbf{5}\\
1000 & 20 & $p_1s_1$ & 1570.94 & 5675.20 & 0.04 & 1 & \textbf{3.16} & \textbf{5674.80} & \textbf{0.00} & \textbf{5}\\
1000 & 20 & $p_1s_2$ & - & 3193.60 & 1.75 & 0 & \textbf{23.58} & \textbf{3148.60} & \textbf{0.00} & \textbf{5}\\
1000 & 20 & $p_1s_3$ & 25.31 & \textbf{7693.80} & \textbf{0.00} & \textbf{5} & \textbf{0.09} & \textbf{7693.80} & \textbf{0.00} & \textbf{5}\\
5000 & 20 & $p_1s_1$ & - & 52548.60 & 522.43 & 0 & \textbf{1.56} & \textbf{28037.80} & \textbf{0.00} & \textbf{5}\\
5000 & 20 & $p_1s_2$ & - & 52520.40 & 100.00 & 0 & \textbf{31.32} & \textbf{15735.40} & \textbf{0.00} & \textbf{5}\\
5000 & 20 & $p_1s_3$ & - & 41035.60 & 6400.14 & 0 & \textbf{0.07} & \textbf{38108.60} & \textbf{0.00} & \textbf{5}\\
\midrule
\multicolumn{11}{c}{Instances with $p_2 = [1,n_J]$ and $B = 20$}\\
\midrule
10 & 20 & $p_2s_1$ & \textbf{0.01} & \textbf{36.20} & \textbf{0.00} & \textbf{5} & 0.05 & \textbf{36.20} & \textbf{0.00} & \textbf{5}\\
10 & 20 & $p_2s_2$ & \textbf{0.02} & \textbf{24.20} & \textbf{0.00} & \textbf{5} & 0.07 & \textbf{24.20} & \textbf{0.00} & \textbf{5}\\
10 & 20 & $p_2s_3$ & \textbf{0.00} & \textbf{42.00} & \textbf{0.00} & \textbf{5} & 0.00 & \textbf{42.00} & \textbf{0.00} & \textbf{5}\\
50 & 20 & $p_2s_1$ & \textbf{0.14} & \textbf{690.20} & \textbf{0.00} & \textbf{5} & 2.52 & \textbf{690.20} & \textbf{0.00} & \textbf{5}\\
50 & 20 & $p_2s_2$ & \textbf{47.90} & \textbf{423.40} & \textbf{0.00} & \textbf{5} & 110.56 & \textbf{423.40} & \textbf{0.00} & \textbf{5}\\
50 & 20 & $p_2s_3$ & \textbf{0.02} & \textbf{1044.60} & \textbf{0.00} & \textbf{5} & 0.06 & \textbf{1044.60} & \textbf{0.00} & \textbf{5}\\
100 & 20 & $p_2s_1$ & \textbf{1.82} & \textbf{2849.40} & \textbf{0.00} & \textbf{5} & 24.15 & \textbf{2849.40} & \textbf{0.00} & \textbf{5}\\
100 & 20 & $p_2s_2$ & - & \textbf{1617.60} & \textbf{0.52} & \textbf{0} & - & \textbf{1617.60} & 0.64 & \textbf{0}\\
100 & 20 & $p_2s_3$ & \textbf{0.08} & \textbf{3889.40} & \textbf{0.00} & \textbf{5} & 0.13 & \textbf{3889.40} & \textbf{0.00} & \textbf{5}\\
500 & 20 & $p_2s_1$ & \textbf{504.40} & \textbf{69070.20} & \textbf{0.02} & \textbf{4} & 994.67 & 69071.00 & 0.02 & 3\\
500 & 20 & $p_2s_2$ & - & 38695.00 & 1.50 & \textbf{0} & - & \textbf{38433.20} & \textbf{0.86} & \textbf{0}\\
500 & 20 & $p_2s_3$ & \textbf{6.18} & \textbf{91883.00} & \textbf{0.00} & \textbf{5} & 30.57 & \textbf{91883.00} & \textbf{0.00} & \textbf{5}\\
1000 & 20 & $p_2s_1$ & - & 272294.40 & 0.06 & 0 & \textbf{1423.96} & \textbf{272187.80} & \textbf{0.01} & \textbf{3}\\
1000 & 20 & $p_2s_2$ & - & 155075.60 & 2.58 & \textbf{0} & - & \textbf{153297.40} & \textbf{1.46} & \textbf{0}\\
1000 & 20 & $p_2s_3$ & \textbf{32.09} & \textbf{371483.00} & \textbf{0.00} & \textbf{5} & 65.80 & \textbf{371483.00} & \textbf{0.00} & \textbf{5}\\
5000 & 20 & $p_2s_1$ & - & 12503592.80 & 421.05 & \textbf{0} & - & \textbf{6659163.00} & \textbf{0.16} & \textbf{0}\\
5000 & 20 & $p_2s_2$ & - & 12495060.40 & 100.00 & \textbf{0} & - & \textbf{4175411.50} & \textbf{9.47} & \textbf{0}\\
5000 & 20 & $p_2s_3$ & - & \textbf{8993880.00} & \textbf{0.00} & \textbf{0} & - & 8994409.60 & 0.01 & \textbf{0}\\
\toprule
\end{tabular}
}
\label{Table_Results1_New_1}
\end{table}

\begin{table}[!pt] 
\centering
\scalefont{1}
\setlength{\tabcolsep}{2.3pt}
\caption{Computational results for new instances proposed ($B=50$).}{
\begin{tabular}{ c c c | r r r r | r r r r }
\toprule
\multicolumn{3}{c}{Instance} & \multicolumn{4}{c}{(MILP$_1^+$)} & \multicolumn{4}{c}{(FLOW)} \\
\cmidrule(l){4-7} \cmidrule(l){8-11}
\multicolumn{1}{c}{Jobs} & 
\multicolumn{1}{c}{Capacity} & 
\multicolumn{1}{c}{Type} &
\multicolumn{1}{c}{$T(s)$} &
\multicolumn{1}{c}{$C_{\max}$} &
\multicolumn{1}{c}{Gap} & 
\multicolumn{1}{c}{\#O} &
\multicolumn{1}{c}{$T(s)$} & 
\multicolumn{1}{c}{$C_{\max}$} & 
\multicolumn{1}{c}{Gap} &
\multicolumn{1}{c}{\#O} \\ 
\midrule
\multicolumn{11}{c}{Instances with $p_1 = [1,20]$ and $B = 50$}\\
\midrule
10 & 50 & $p_1s_1$ & \textbf{0.01} & \textbf{65.80} & \textbf{0.00} & \textbf{5} & 0.13 & \textbf{65.80} & \textbf{0.00} & \textbf{5}\\
10 & 50 & $p_1s_2$ & \textbf{0.02} & \textbf{40.60} & \textbf{0.00} & \textbf{5} & 0.18 & \textbf{40.60} & \textbf{0.00} & \textbf{5}\\
10 & 50 & $p_1s_3$ & \textbf{0.00} & \textbf{109.20} & \textbf{0.00} & \textbf{5} & 0.00 & \textbf{109.20} & \textbf{0.00} & \textbf{5}\\
50 & 50 & $p_1s_1$ & \textbf{0.22} & \textbf{299.80} & \textbf{0.00} & \textbf{5} & 7.79 & \textbf{299.80} & \textbf{0.00} & \textbf{5}\\
50 & 50 & $p_1s_2$ & 72.71 & \textbf{169.00} & \textbf{0.00} & \textbf{5} & \textbf{53.15} & \textbf{169.00} & \textbf{0.00} & \textbf{5}\\
50 & 50 & $p_1s_3$ & \textbf{0.01} & \textbf{429.00} & \textbf{0.00} & \textbf{5} & 0.05 & \textbf{429.00} & \textbf{0.00} & \textbf{5}\\
100 & 50 & $p_1s_1$ & \textbf{1.01} & \textbf{560.40} & \textbf{0.00} & \textbf{5} & 11.33 & \textbf{560.40} & \textbf{0.00} & \textbf{5}\\
100 & 50 & $p_1s_2$ & - & \textbf{317.80} & 0.93 & 0 & \textbf{654.39} & \textbf{317.80} & \textbf{0.08} & \textbf{4}\\
100 & 50 & $p_1s_3$ & \textbf{0.05} & \textbf{746.60} & \textbf{0.00} & \textbf{5} & 0.08 & \textbf{746.60} & \textbf{0.00} & \textbf{5}\\
500 & 50 & $p_1s_1$ & 1465.28 & 2741.20 & 0.06 & 2 & \textbf{84.02} & \textbf{2740.60} & \textbf{0.00} & \textbf{5}\\
500 & 50 & $p_1s_2$ & - & 1602.40 & 2.77 & 0 & \textbf{1625.40} & \textbf{1572.80} & \textbf{0.25} & \textbf{1}\\
500 & 50 & $p_1s_3$ & 1.24 & \textbf{4036.60} & \textbf{0.00} & \textbf{5} & \textbf{0.25} & \textbf{4036.60} & \textbf{0.00} & \textbf{5}\\
1000 & 50 & $p_1s_1$ & - & 5519.80 & 26.27 & 0 & \textbf{34.75} & \textbf{5517.00} & \textbf{0.00} & \textbf{5}\\
1000 & 50 & $p_1s_2$ & - & 3219.60 & 3.80 & 0 & \textbf{1654.71} & \textbf{3113.60} & \textbf{0.21} & \textbf{1}\\
1000 & 50 & $p_1s_3$ & 9.58 & \textbf{7901.00} & \textbf{0.00} & \textbf{5} & \textbf{0.14} & \textbf{7901.00} & \textbf{0.00} & \textbf{5}\\
5000 & 50 & $p_1s_1$ & - & 52744.20 & 186.78 & 0 & \textbf{57.91} & \textbf{27336.00} & \textbf{0.00} & \textbf{5}\\
5000 & 50 & $p_1s_2$ & - & 52553.40 & 100.00 & 0 & \textbf{1581.85} & \textbf{15769.60} & \textbf{0.03} & \textbf{1}\\
5000 & 50 & $p_1s_3$ & 1766.61 & 38921.60 & 3400.02 & 1 & \textbf{0.24} & \textbf{38914.80} & \textbf{0.00} & \textbf{5}\\
\midrule
\multicolumn{11}{c}{Instances with $p_2 = [1,n_J]$ and $B = 50$}\\
\midrule
10 & 50 & $p_2s_1$ & \textbf{0.00} & \textbf{46.80} & \textbf{0.00} & \textbf{5} & 0.04 & \textbf{46.80} & \textbf{0.00} & \textbf{5}\\
10 & 50 & $p_2s_2$ & \textbf{0.04} & \textbf{23.20} & \textbf{0.00} & \textbf{5} & 0.13 & \textbf{23.20} & \textbf{0.00} & \textbf{5}\\
10 & 50 & $p_2s_3$ & \textbf{0.00} & \textbf{48.80} & \textbf{0.00} & \textbf{5} & 0.00 & \textbf{48.80} & \textbf{0.00} & \textbf{5}\\
50 & 50 & $p_2s_1$ & \textbf{0.22} & \textbf{780.60} & \textbf{0.00} & \textbf{5} & 20.29 & \textbf{780.60} & \textbf{0.00} & \textbf{5}\\
50 & 50 & $p_2s_2$ & \textbf{89.05} & \textbf{411.40} & \textbf{0.00} & \textbf{5} & 740.47 & \textbf{411.40} & \textbf{0.00} & \textbf{5}\\
50 & 50 & $p_2s_3$ & \textbf{0.01} & \textbf{1010.40} & \textbf{0.00} & \textbf{5} & 0.09 & \textbf{1010.40} & \textbf{0.00} & \textbf{5}\\
100 & 50 & $p_2s_1$ & \textbf{4.58} & \textbf{2640.00} & \textbf{0.00} & \textbf{5} & 1028.85 & 2640.20 & 0.10 & 3 \\
100 & 50 & $p_2s_2$ & \textbf{1698.37} & 1636.00 & \textbf{0.68} & \textbf{1} & - & \textbf{1634.00} & 1.06 & 0\\
100 & 50 & $p_2s_3$ & \textbf{0.05} & \textbf{4251.00} & \textbf{0.00} & \textbf{5} & 0.21 & \textbf{4251.00} & \textbf{0.00} & \textbf{5}\\
500 & 50 & $p_2s_1$ & \textbf{899.33} & \textbf{65347.80} & \textbf{0.03} & \textbf{3} & - & 65388.80 & 0.13 & 0 \\
500 & 50 & $p_2s_2$ & - & 39222.40 & 2.71 & \textbf{0} & - & \textbf{38819.00} & \textbf{1.37} & \textbf{0}\\
500 & 50 & $p_2s_3$ & \textbf{1.23} & \textbf{94798.00} & \textbf{0.00} & \textbf{5} & 10.94 & \textbf{94798.00} & \textbf{0.00} & \textbf{5}\\
1000 & 50 & $p_2s_1$ & \textbf{1555.41} & 262289.20 & 0.26 & \textbf{1} & 1711.02 & \textbf{261181.50} & \textbf{0.16} & \textbf{1}\\
1000 & 50 & $p_2s_2$ & - & \textbf{157692.20} & \textbf{4.07} & \textbf{0} & - & No solution &\multicolumn{1}{c}{$\infty$} & \textbf{0}\\
1000 & 50 & $p_2s_3$ & \textbf{7.95} & \textbf{371115.80} & \textbf{0.00} & \textbf{5} & 48.52 & \textbf{371115.80} & \textbf{0.00} & \textbf{5}\\
5000 & 50 & $p_2s_1$ & - & \textbf{12573122.60} & \textbf{161.91} & \textbf{0} & - & No solution &\multicolumn{1}{c}{$\infty$} & \textbf{0}\\
5000 & 50 & $p_2s_2$ & - & \textbf{12546134.00} & \textbf{100.00} & \textbf{0} & - & No solution &\multicolumn{1}{c}{$\infty$} & \textbf{0}\\
5000 & 50 & $p_2s_3$ & 797.53 & \textbf{9291848.20} & \textbf{0.00} & \textbf{4} & \textbf{524.10} & 9291858.80 & \textbf{0.00} & \textbf{4}\\
\toprule
\end{tabular}
}
\label{Table_Results1_New_2}
\end{table}

\begin{table}[!pt] 
\centering
\scalefont{1}
\setlength{\tabcolsep}{2.3pt}
\caption{Computational results for  new instances proposed ($B=100$).}{
\begin{tabular}{ c c c | r r r r | r r r r }
\toprule
\multicolumn{3}{c}{Instance} & \multicolumn{4}{c}{(MILP$_1^+$)} & \multicolumn{4}{c}{(FLOW)} \\
\cmidrule(l){4-7} \cmidrule(l){8-11}
\multicolumn{1}{c}{Jobs} & 
\multicolumn{1}{c}{Capacity} & 
\multicolumn{1}{c}{Type} &
\multicolumn{1}{c}{$T(s)$} &
\multicolumn{1}{c}{$C_{\max}$} &
\multicolumn{1}{c}{Gap} & 
\multicolumn{1}{c}{\#O} &
\multicolumn{1}{c}{$T(s)$} & 
\multicolumn{1}{c}{$C_{\max}$} & 
\multicolumn{1}{c}{Gap} &
\multicolumn{1}{c}{\#O} \\ 
\midrule
\multicolumn{11}{c}{Instances with $p_1 = [1,20]$ and $B = 100$}\\
\midrule
10 & 100 & $p_1s_1$ & \textbf{0.01} & \textbf{77.20} & \textbf{0.00} & \textbf{5} & 0.06 & \textbf{77.20} & \textbf{0.00} & \textbf{5}\\
10 & 100 & $p_1s_2$ & \textbf{0.01} & \textbf{44.00} & \textbf{0.00} & \textbf{5} & 0.24 & \textbf{44.00} & \textbf{0.00} & \textbf{5}\\
10 & 100 & $p_1s_3$ & \textbf{0.00} & \textbf{94.80} & \textbf{0.00} & \textbf{5} & 0.01 & \textbf{94.80} & \textbf{0.00} & \textbf{5}\\
50 & 100 & $p_1s_1$ & \textbf{0.07} & \textbf{316.00} & \textbf{0.00} & \textbf{5} & 3.39 & \textbf{316.00} & \textbf{0.00} & \textbf{5}\\
50 & 100 & $p_1s_2$ & \textbf{9.49} & \textbf{180.60} & \textbf{0.00} & \textbf{5} & 189.28 & \textbf{180.60} & \textbf{0.00} & \textbf{5}\\
50 & 100 & $p_1s_3$ & \textbf{0.01} & \textbf{411.20} & \textbf{0.00} & \textbf{5} & 0.07 & \textbf{411.20} & \textbf{0.00} & \textbf{5}\\
100 & 100 & $p_1s_1$ & \textbf{2.05} & \textbf{586.80} & \textbf{0.00} & \textbf{5} & 82.86 & \textbf{586.80} & \textbf{0.00} & \textbf{5}\\
100 & 100 & $p_1s_2$ & \textbf{1372.59} & 329.00 & 0.74 & \textbf{2} & 1588.41 & \textbf{328.60} & \textbf{0.47} & \textbf{2}\\
100 & 100 & $p_1s_3$ & \textbf{0.05} & \textbf{811.20} & \textbf{0.00} & \textbf{5} & 0.16 & \textbf{811.20} & \textbf{0.00} & \textbf{5}\\
500 & 100 & $p_1s_1$ & 877.55 & 2715.40 & 0.07 & 3 & \textbf{664.78} & \textbf{2714.40} & \textbf{0.01} & \textbf{4}\\
500 & 100 & $p_1s_2$ & - & 1601.60 & 3.09 & \textbf{0} & - & \textbf{1595.20} & \textbf{2.26} & \textbf{0}\\
500 & 100 & $p_1s_3$ & 1.71 & \textbf{4096.20} & \textbf{0.00} & \textbf{5} & \textbf{1.06} & \textbf{4096.20} & \textbf{0.00} & \textbf{5}\\
1000 & 100 & $p_1s_1$ & - & 5591.40 & 0.63 & 0 & \textbf{902.17} & \textbf{5569.40} & \textbf{0.02} & \textbf{3}\\
1000 & 100 & $p_1s_2$ & - & 3282.00 & 4.51 & \textbf{0} & - & \textbf{3196.80} & \textbf{1.72} & \textbf{0}\\
1000 & 100 & $p_1s_3$ & 13.43 & \textbf{8094.80} & \textbf{0.00} & \textbf{5} & \textbf{1.46} & \textbf{8094.80} & \textbf{0.00} & \textbf{5}\\
5000 & 100 & $p_1s_1$ & - & 52416.60 & 123.69 & 0 & \textbf{883.92} & \textbf{26898.60} & \textbf{0.00} & \textbf{3}\\
5000 & 100 & $p_1s_2$ & - & 52518.60 & 100.00 & \textbf{0} & - & \textbf{15879.40} & \textbf{0.62} & \textbf{0}\\
5000 & 100 & $p_1s_3$ & 1481.01 & \textbf{39307.40} & \textbf{0.00} & \textbf{5} & \textbf{0.62} & \textbf{39307.40} & \textbf{0.00} & \textbf{5}\\
\midrule
\multicolumn{11}{c}{Instances with $p_2 = [1,n_J]$ and $B = 100$}\\
\midrule
10 & 100 & $p_2s_1$ & \textbf{0.00} & \textbf{32.60} & \textbf{0.00} & \textbf{5} & 0.10 & \textbf{32.60} & \textbf{0.00} & \textbf{5}\\
10 & 100 & $p_2s_2$ & \textbf{0.02} & \textbf{21.40} & \textbf{0.00} & \textbf{5} & 0.18 & \textbf{21.40} & \textbf{0.00} & \textbf{5}\\
10 & 100 & $p_2s_3$ & \textbf{0.00} & \textbf{49.40} & \textbf{0.00} & \textbf{5} & 0.00 & \textbf{49.40} & \textbf{0.00} & \textbf{5}\\
50 & 100 & $p_2s_1$ & \textbf{0.16} & \textbf{823.60} & \textbf{0.00} & \textbf{5} & 281.80 & \textbf{823.60} & \textbf{0.00} & \textbf{5}\\
50 & 100 & $p_2s_2$ & \textbf{40.72} & \textbf{441.60} & \textbf{0.00} & \textbf{5} & 778.21 & \textbf{441.60} & 0.40 & 4\\
50 & 100 & $p_2s_3$ & \textbf{0.01} & \textbf{1041.40} & \textbf{0.00} & \textbf{5} & 0.13 & \textbf{1041.40} & \textbf{0.00} & \textbf{5}\\
100 & 100 & $p_2s_1$ & \textbf{0.63} & \textbf{3092.20} & \textbf{0.00} & \textbf{5} & 654.44 & 3092.40 & 0.02 & 4\\
100 & 100 & $p_2s_2$ & - & \textbf{1631.60} & \textbf{0.59} & \textbf{0} & - & 1636.40 & 1.59 & \textbf{0}\\
100 & 100 & $p_2s_3$ & \textbf{0.06} & \textbf{4090.20} & \textbf{0.00} & \textbf{5} & 0.68 & \textbf{4090.20} & \textbf{0.00} & \textbf{5}\\
500 & 100 & $p_2s_1$ & \textbf{1207.01} & \textbf{65728.00} & \textbf{0.02} & \textbf{2} & - & 69647.00 & 1.32 & 0\\
500 & 100 & $p_2s_2$ & - & \textbf{39355.60} & \textbf{3.43} & \textbf{0} & - & 41833.00 & 9.15 & \textbf{0}\\
500 & 100 & $p_2s_3$ & \textbf{0.76} & \textbf{98099.00} & \textbf{0.00} & \textbf{5} & 6.06 & \textbf{98099.00} & \textbf{0.00} & \textbf{5}\\
1000 & 100 & $p_2s_1$ & - & \textbf{262033.20} & \textbf{0.32} & \textbf{0} & - & No solution &\multicolumn{1}{c}{$\infty$} & \textbf{0}\\
1000 & 100 & $p_2s_2$ & - & \textbf{168105.40} & \textbf{9.53} & \textbf{0} & - & 178472.20 & 14.76 & \textbf{0}\\
1000 & 100 & $p_2s_3$ & \textbf{6.72} & \textbf{381506.40} & \textbf{0.00} & \textbf{5} & 46.58 & \textbf{381506.40} & \textbf{0.00} & \textbf{5}\\
5000 & 100 & $p_2s_1$ & - & \textbf{12529354.80} & \textbf{107.49} & \textbf{0} & - & No solution &\multicolumn{1}{c}{$\infty$} & \textbf{0}\\
5000 & 100 & $p_2s_2$ & - & \textbf{12488031.00} & \textbf{100.00} & \textbf{0} & - & No solution &\multicolumn{1}{c}{$\infty$} & \textbf{0}\\
5000 & 100 & $p_2s_3$ & \textbf{891.11} & \textbf{9404104.80} & \textbf{0.00} & \textbf{4} & 1493.59 & 9404178.20 & \textbf{0.00} & 1\\
\toprule
\end{tabular}
}
\label{Table_Results1_New_3}
\end{table}

The comparative tests show that new instances are more difficult for (FLOW), especially when the  processing times are selected from interval  $p_2$, as there are more distinct processing times in this case. In this way, (FLOW) needs to consider many flow structures, one for each processing time.
The difference between the instances $p_1$ and $p_2$, however, are not influencing the number of variables in (MILP$_1^+$), and the computational results show that this model is more stable when the range of the processing time changes.

(MILP$_1^+$) is superior for instances of type $p_2$ especially when the machine capacity  increases. 
The results in Table \ref{Table_Results1_New_2} show that (FLOW) was not able to find even a single integer solution in some cases. 
This behavior occurs in the worst case scenario for (FLOW), where a large value for the machine capacity is combined with the type of instances $p_2$. When $B = 100$, the (FLOW) model considers many nodes in the arc-flow structures, which reflects in the computational performance.
Another difference between this new set of instances and the instances tested in the previous section is that the range of job sizes varies with the number of jobs. The number of arcs in (FLOW) increases as the number of jobs increases.

Even with the difficulties created with this new set of instances, the (FLOW) model obtained superior results for instances of type $p_1$, especially when the number of jobs increases. 
Considering the instances of type $p_1$ with $B$ up to 20, it was possible to find the optimal solution for all instances.
These tests show that both models are complementary, that is, there are clearly situations in which each model is superior to the other.

\section{Conclusions and Future Work} \label{chap:7}
\label{sec:conc}
We address a single batch processing machine scheduling problem. The economic importance of the problem has motivated the investigation of good solution approaches, and its NP-hardness has led the majority of this research to focus on heuristic approaches. 
We show that applying a good \ac{MILP} formulation for this scheduling problem we can go a step further in the exact resolution of applied problems, having presented optimal solutions for test instances with sizes never considered in the literature by exact methods. We were able to present much better  results for a set of benchmark instances  from the literature, with the formulation proposed, when compared to results obtained  with two other formulations from the literature.

We propose an arc flow-based formulation for the  problem, in which the numbers of variables and constraints do not change when the number of jobs increases. This procedure enabled us to find the optimum solution of an instance with 100 million jobs in 3.25 seconds, and define a new threshold for the size of the instances, as the maximum number of jobs previously treated in the literature was 500, using heuristic approaches for instances with the same parameter settings. 

Analyzing the parameters that affect the number of variables and constraints in our arc-flow formulation, we conclude that the formulation  is very effective for instances of the problem, where, although the number of jobs is very large, they have similar characteristics, so the numbers of distinct processing times and job sizes are not very large.  In order to investigate the scability of our approach, we have generated a new set of instances, where the parameters setting had the purpose of making the instances more challenging for the model than the benchmark instances from the literature. We see that the our model complements another model from the literature, in the sense that depending on the characteristics of the instances, each one can perform better than the other. These characteristics are identified and can be found on different applications for the problem. For example, in burn-in tests for semiconductors, it is natural to have many jobs to perform with similar characteristics of processing times and sizes.  These tests were mentioned in the literature as an important   application for the single batch processing machine scheduling problem.

As future research, we would like to investigate if the good performance of the formulation  presented can be replicated when the arc flow-based approach is applied to other problems in the vast area of scheduling applications as, for example, the problem of scheduling a batch processing machine with incompatible job families.

\bibliographystyle{plain}
\bibliography{references}
\end{document}